\documentclass[titlepage,11pt]{article}
\usepackage{latexsym}
\usepackage{amsfonts}
\usepackage{amsmath}
\usepackage{enumitem}
\usepackage{tikz}
\usepackage{float}
\usepackage{lmodern}
\usepackage{tikz}
\usetikzlibrary{decorations.pathmorphing}
\usetikzlibrary{graphs}

\newcommand\blackslug{\hbox{\hskip 1pt \vrule width 4pt height 8pt depth 1.5pt
        \hskip 1pt}}
\newcommand\bbox{\hfill \quad \blackslug \bigbreak}

\renewcommand{\epsilon}{\varepsilon}
\def\cc{\hbox{-}\cdots\hbox{-}}

\def\ll{,\ldots,}

\title{Induced subgraphs of graphs with large chromatic number.
\\XIII. New brooms}
\author{Alex Scott\thanks{Supported by a Leverhulme Research
Fellowship}\\
Mathematical Institute, University of Oxford, Oxford OX2 6GG, UK
\\
\\
Paul Seymour\thanks{Supported by ONR grant N00014-14-1-0084 and NSF
grant DMS-1265563.}\\
Princeton University, Princeton, NJ 08544, USA}

\date{October 24, 2016; revised \today}

\newtheorem{thm}{}[section]

\newcommand{\Proof}{\noindent{\bf Proof.}\ \ }

\begin{document}
\maketitle
\begin{abstract}
Gy\'arf\'as~\cite{gyarfas} and Sumner~\cite{sumner} independently conjectured that for every tree $T$,
the class of graphs not containing $T$ as an induced subgraph is $\chi$-bounded, that is, the chromatic numbers of graphs in this class
are bounded above by a function of their clique numbers. This remains open for general trees $T$, but has been proved for some
particular trees.
For $k\ge 1$, let us say a {\em broom}
of {\em length $k$}
is a tree obtained from a $k$-edge path with ends $a,b$ by adding some number of leaves adjacent to $b$, and we call $a$
its {\em handle}. A tree obtained from brooms of lengths $k_1\ll k_n$ by identifying their handles is a {\em $(k_1\ll k_n)$-multibroom}.
Kierstead and Penrice~\cite{kierstead1} proved that every $(1\ll 1)$-multibroom $T$ satisfies the Gy\'arf\'as-Sumner conjecture, 
and Kierstead and Zhu~\cite{kierstead2} proved the same for
$(2\ll 2)$-multibrooms. 

In this paper give a common generalization; we prove that every $(1\ll 1,2\ll 2)$-multibroom satisfies the Gy\'arf\'as-Sumner 
conjecture .

\end{abstract}

\section{Introduction}
For a graph $G$, let $\chi(G)$ denote the chromatic number of $G$, and let $\omega(G)$ denote its clique number, that is,  
the number of vertices in its largest clique. We say a graph $G$ {\em contains} $H$ if some induced subgraph of $G$ is isomorphic to $H$,
and otherwise $G$ is {\em $H$-free}.

The Gy\'arf\'as-Sumner conjecture~\cite{gyarfas,sumner} asserts that:
\begin{thm}\label{conj}{\bf Conjecture: } For every forest $T$ and every integer $\kappa$, there exists $c$ such that $\chi(G)\le c$ for
every $T$-free graph $G$ with $\omega(G)\le \kappa$.
\end{thm}

There has been surprisingly little progress on this conjecture. It is easy to see that if the conjecture holds for every component 
of a forest then it holds for the forest (the first component must be present; delete it and all vertices with a neighbour in it and 
repeat with the next component), and so it suffices to prove the conjecture when $T$ is a tree. Gy\'arf\'as~\cite{gyarfas}
proved the conjecture when $T$ is a path, and 
Scott~\cite{scott} proved it when $T$ is a subdivision of a star; and recently, 
with Maria Chudnovsky, we~\cite{cuttingedge} proved it for trees obtained
from a subdivided star by adding one more vertex with one neighbour, and for trees obtained from a star and a subdivided star
by adding a path between their centres. But the results that concern us most here are theorems of
Gy\'arf\'as, Szemer\'edi and Tuza~\cite{gyarfas2}, Kierstead and Penrice~\cite{kierstead1}, and 
Kierstead and Zhu~\cite{kierstead2}, which are the only other results so far on the Gy\'arf\'as-Sumner conjecture, and which we explain
next.

For $k\ge 1$, let us say a {\em broom}
of {\em length $k$}
is a tree obtained from a $k$-edge path with ends $a,b$ by adding some number of leaves adjacent to $b$, and we call $a$
its {\em handle}. A tree obtained from $n$ brooms of lengths $k_1\ll k_n$ respectively by identifying their handles is called
a {\em $(k_1\ll k_n)$-multibroom}. Gy\'arf\'as, Szemer\'edi and Tuza (in the triangle-free case) and then
Kierstead and Penrice (in the general case) proved that $(1\ll 1)$-multibrooms satisfy the Gy\'arf\'as-Sumner conjecture, and 
Kierstead and Zhu proved that $(2\ll 2)$-multibrooms satisfy it.
In this paper we prove a common generalization of these results:
every $(1\ll 1,2\ll 2)$-multibroom satisfies the Gy\'arf\'as-Sumner conjecture.

\begin{figure}[H]
\centering

\begin{tikzpicture}[scale=.8,auto=left]
\tikzstyle{every node}=[inner sep=1.5pt, fill=black,circle,draw]

\node (x) at (-4.5,-4) {};
\node (i) at (-1.5,-6) {};
\node (j) at ( -3.5,-6) {};
\node (k) at ( -5.5,-6) {};
\node (l) at ( -7.5,-6) {};
\node (i1) at (-0.9, -7) {};
\node (i2) at (-1.5, -7) {};
\node (i3) at (-2.1, -7) {};
\node (j1) at (-2.9, -7) {};
\node (j2) at (-3.5, -7) {};
\node (j3) at (-4.1, -7) {};
\node (k1) at (-4.9, -7) {};
\node (k2) at (-5.5, -7) {};
\node (k3) at (-6.1, -7) {};
\node (l1) at (-6.9, -7) {};
\node (l2) at (-7.5, -7) {};
\node (l3) at (-8.1, -7) {};
\node (ii) at (-3,-5) {};
\node (jj) at (-4,-5) {};

\foreach \from/\to in {x/ii,x/jj,x/k,x/l,ii/i,jj/j,i/i1,i/i2,i/i3,j/j1,j/j2,j/j3,k/k1,k/k2,k/k3,l/l1,l/l2,l/l3}
\draw [-] (\from) -- (\to);

\end{tikzpicture}

\caption{A $(1,1,2,2)$-multibroom} \label{fig:1}
\end{figure}
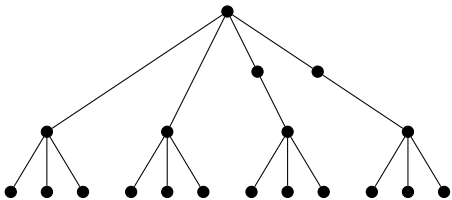

Let us state this
more precisely.
A {\em $(k,\delta)$-broom} means a broom of length $k$ with $\delta$ leaves different from its handle (thus, it is obtained by adding
$\delta$ leaves adjacent to one end of a $k$-edge path).
For $\delta\ge 1$, let $T(\delta)$ be the tree formed from the disjoint union of $\delta$
$(1,\delta)$-brooms and $\delta$ $(2,\delta)$-brooms by identifying their handles. We will prove that
\begin{thm}\label{mainthm1}
For all $\delta\ge 0$ and all $\kappa\ge 0$ there exists $c$ such that every $T(\delta)$-free graph with $\omega(G)\le \kappa$
has chromatic number at most $c$.
\end{thm}
The proof method is by combining ideas of~\cite{kierstead1,kierstead2} with some new twists.

\section{Inductions}

There are various inductions that will give us some assistance. We can use induction on $\kappa$, and on $\delta$ (in fact with a little
work we can more-or-less assume that the result holds for every tree obtained from $T(\delta)$ by deleting a leaf), and there is a third induction,
core maximization, that we explain later. Next we explain these inductions in more detail.

First and easiest, by induction on $\kappa$, we may assume that there exists $\tau$ such that $\chi(G)\le \tau$
for every $T(\delta)$-free graph with clique number less than $\kappa$. In particular, this tells us that
if $G$ is $T(\delta)$-free with clique number at most $\delta$, then for every vertex, the subgraph induced
on its neighbours has chromatic number at most $\tau$ (since this subgraph has clique number less than $\tau$).
Consequently we can use \ref{rad2} below, taking $T=T(\delta)$.

If $v$ is a vertex of a graph $G$, and $k\ge 1$, $N^{k}(v)$ or $N^{k}_G(v)$ denotes the set of vertices of $G$ 
with distance exactly $k$
from $v$, and $N^{k}[v]$ denotes the set with distance at most $k$ from $v$.
If $G$ is a nonnull graph  and $k\ge 1$,
we define $\chi^{k}(G)$ to be the maximum of $\chi(N^{k}[v])$ taken over all vertices $v$ of $G$.
(For the null graph $G$ we define $\chi^{k}(G)=0$.) 

The following follows by repeated application of theorem 3.2 of~\cite{cuttingedge} (a similar theorem for $(2\ll 2)$-multibrooms
is proved in \cite{kierstead2}):
\begin{thm}\label{rad2}
Let $T$ be a tree formed by identifying the handles of some set of brooms (of arbitrary lengths).
For all $\kappa, \tau\ge 0$ there exists $c$
with the following property.  Let $G$ be a $T$-free graph, with $\omega(G)\le \kappa$,
such that for every vertex, the subgraph induced
on its neighbours has chromatic number at most $\tau$.
Then $\chi^2(G)\le c$.
\end{thm}

Next, let us explore induction on the size of $T(\delta)$. That will allow us to exploit ``matching-covered'' sets.
Let $X\subseteq V(G)$. 
We say that $X$ is {\em matching-covered} in $G$ if for each $x\in X$ there exists $y\in V(G)\setminus X$
adjacent to $x$ and to no other vertex in $X$.

We would like to be able to assume that the result holds for all trees obtained from $T(\delta)$ by deleting a leaf;
but only deleting one leaf, from one of its brooms, and so the smaller tree is not of the form $T(\delta')$ for $\delta'<\delta$,
and so induction on $\delta$ is not fine enough. We could change the statement of the theorem, and prove it not only for
$T(\delta)$, but for any tree that is a subtree of $T(\delta)$; but that would make things notationally more complicated later.  
There is another way to do it that is more convenient.

Let us say that $G$ is {\em $(\delta,\kappa)$-good} if
$G$ is $T(\delta)$-free and $\omega(G)\le \kappa$. An {\em ideal} of graphs is a class $\mathcal{C}$ of graphs
such that every induced subgraph of a member of $\mathcal{C}$ also belongs to $\mathcal{C}$.
If $X\subseteq V(G)$, we write $\chi(X)$ for $\chi(G[X])$ when there is no ambiguity.

\begin{thm}\label{matchingcovered}
Let $\delta,\kappa\ge 1$. If $\mathcal{C}$ is an ideal of $(\delta,\kappa)$-good graphs with unbounded chromatic number, then there
exist a subideal $\mathcal{C}'$ of $\mathcal{C}$  with unbounded chromatic number, and a number $c$ such that 
every matching-covered set in every member of $\mathcal{C}'$ has chromatic number at most $c$.
\end{thm}
\Proof
Let $R$ be a maximal subtree of $T(\delta)$ such that
there exists $c$ such that 
every $R$-free graph in $\mathcal{C}$ with clique number at most $\kappa$ has chromatic number at most $c$, 
and choose some such number $c$.
By hypothesis there are members of $\mathcal{C}$ with arbitrarily large chromatic number, so 
$R\ne T(\delta)$.
Hence there is a subtree $S$ of $T(\delta)$ with a leaf $v$
such that $S\setminus v = R$. Let $u$ be the neighbour of $v$ in $S$.
Let $\mathcal{C}'$ be the subideal of all $S$-free graphs in $\mathcal{C}$.
From the maximality of $R$, there are graphs in $\mathcal{C}$ with arbitrarily large chromatic number. 

Let $G\in \mathcal{C}$, and let $X$ be matching-covered in $G$. Suppose that
there is an induced subgraph of $G[X]$ isomorphic to $R$, and to simplify notation we assume
it equals $R$. Choose $y\in V(G)\setminus X$ adjacent to $u$ and to no other vertex in $X$; then $G[V(R)\cup \{y\}]$
is isomorphic to $S$, a contradiction. Thus $G[X]$ does not contain $R$. Since $G[X]\in \mathcal{C}$,
the choice of $c$ implies 
that $\chi(X)\le c$. Since all graphs in $\mathcal{C}$ are $(\delta,\kappa)$-good, this proves \ref{matchingcovered}.~\bbox

There is a third, very helpful, induction we can use, but it is more complicated.
For integers $a,b\ge 1$, let us say an {\em $(a,b)$-core} in a graph $G$ is a subset $Y\subseteq V(G)$
of cardinality $ab$, that admits a partition $\{A_1\ll A_b\}$ 
such that 
\begin{itemize}
\item $A_1\ll A_b$ each have cardinality $a$;
\item $A_1\ll A_b$ are all stable sets of $G$; and
\item for $1\le i<j\le b$, every vertex in $A_i$ is adjacent to every vertex in $A_j$.
\end{itemize}
(An $(a,b)$-core is therefore a complete multipartite induced subgraph of specified size.)
This partition is unique, since $a\ge 1$, and we speak of $A_1\ll A_b$ as the {\em parts} of $Y$.
Thus, if there is an $(a,b)$-core in $G$ then $b\le \omega(G)$. 
Let $\mathbb{N}$ denote the set of nonnegative integers.

\begin{thm}\label{cores}
Let $\delta, \kappa\ge 1$. If $\mathcal{C}$ is an ideal of $(\delta,\kappa)$-good graphs with unbounded chromatic number, then there
exist a subideal $\mathcal{C}'$ of $\mathcal{C}$  with unbounded chromatic number, and integers $\alpha\ge 1$ and
$\beta\ge 2$, and a non-decreasing function $\theta:\mathbb{N}\rightarrow \mathbb{N}$,
 with the following properties:
\begin{itemize}
\item for all $a\ge 1$,
every graph in $\mathcal{C}$ with chromatic number more than $\theta(a)$ admits an $(a,\beta)$-core.
\item no graph in $\mathcal{C}'$ admits an $(\alpha,\beta+1)$-core.
\end{itemize}
\end{thm}
\Proof
For integers $a\ge 1$ and $b\ge 2$, let us say $(a,b)$ is {\em unavoidable} if there exists $c$ such that 
every graph in $\mathcal{C}$ with chromatic number more than $c$ admits an $(a,b)$-core.
V. R\"odl (see~\cite{rodl}) proved that for all integers $a\ge 1$, $(a,2)$ is unavoidable. 
Choose $\beta$ with $2\le \beta \le \kappa+1$ maximum such that for all $a\ge 1$, $(a,\beta)$ is unavoidable.
Since $(a,\beta)$ is unavoidable for all $a\ge 1$, there is a function       
$\theta:\mathbb{N}\rightarrow \mathbb{N}$ such that for all $a\ge 1$, 
every graph in $\mathcal{C}$ with chromatic number more than $\theta(a)$ admits an $(a,\beta)$-core,
and we can choose $\theta$ to be non-decreasing, so the first bullet holds.

By hypothesis there are graphs in $\mathcal{C}$ with unbounded chromatic number, and they do not admit
$(1,\kappa+1)$-cores (because they have clique number at most $\kappa$), so $\beta\le \kappa$. 
From the maximality of $\beta$,
there exists $\alpha\ge 1$
such that there are graphs in $\mathcal{C}$ with arbitrarily large chromatic number that do not admit an $(\alpha,\beta+1)$-core.
Let $\mathcal{C}'$ be the ideal of graphs in $\mathcal{C}$ that do not admit an $(\alpha,\beta+1)$-core; then
the second bullet holds. This proves \ref{cores}.~\bbox

We combine these results in the following.
\begin{thm}\label{induction}
Let $\delta\ge 1$. Suppose that for some value of $\kappa\ge 1$ there are $(\delta,\kappa)$-good graphs with unbounded 
chromatic number. Then there exist $\tau\ge 0$, $\alpha\ge 1$, $\beta,\kappa\ge 2$, 
a non-decreasing function $\theta:\mathbb{N}\rightarrow \mathbb{N}$,
and an ideal $\mathcal{C}$ of graphs with unbounded chromatic number, such that
for every $G\in \mathcal{C}$:
\begin{itemize}
\item $G$ is $T(\delta)$-free;
\item $\omega(G)\le \kappa$;
\item $\chi^2(G)\le \tau$;
\item every matching-covered set in $G$ has chromatic number at most $\tau$;
\item for all $a\ge 1$, if $\chi(G)>\theta(a)$
then $G$ admits an $(a,\beta)$-core;
\item $G$ does not admit an $(\alpha,\beta+1)$-core.
\end{itemize}
\end{thm}
\Proof
Choose $\kappa$ minimum such that there are $(\delta,\kappa)$-good graphs with unbounded 
chromatic number. Thus $\kappa\ge 2$. Choose $\tau_1$ such that every $(\delta,\kappa-1)$-good graph
has chromatic number at most $\tau_1$. By \ref{rad2} there exists $\tau_2$ such that $\chi^2(G)\le \tau_2$
for every $(\delta,\kappa)$-good graph. By \ref{matchingcovered} there exist $\tau_3$ and an ideal $\mathcal{C}_1$ 
of $(\delta,\kappa)$-good graphs with unbounded
chromatic number, such that every matching-covered set in $G$ has chromatic number at most $\tau_3$. By \ref{cores},
there exist a subideal $\mathcal{C}$ of $\mathcal{C}_1$ with unbounded
chromatic number, and $\alpha,\beta$ satisfying the last two bullets. Let $\tau = \max(\tau_1,\tau_2,\tau_3)$; then
all six bullets are satisfied. This proves \ref{induction}.~\bbox

In view of \ref{induction}, in order to prove \ref{mainthm1} it suffices to show the following:
\begin{thm}\label{mainthm2}
For all $\tau\ge 0$, and $\alpha, \delta\ge 1$, and $\beta\ge 2$, and for every non-decreasing function 
$\theta:\mathbb{N}\rightarrow \mathbb{N}$, there exists $c$ such that if $G$ satisfies
{\rm 
\begin{enumerate}[label=(\roman*)]
\item {\em $G$ is $T(\delta)$-free;}
\item {\em $\chi^2(G)\le \tau$;}
\item {\em every matching-covered set in $G$ has chromatic number at most $\tau$;}
\item {\em for all $a\ge 1$, if $\chi(G)>\theta(a)$
then $G$ admits an $(a,\beta)$-core;}
\item {\em $G$ does not admit an $(\alpha,\beta+1)$-core}
\end{enumerate}
{\em then $\chi(G)\le c$.}}
\end{thm}
We could have added another constant $\kappa$ and another condition that $\omega(G)\le \kappa$, but it turns out not to be 
needed any more (a bound on $\omega(G)$ is implied by the second condition).

The five statements (i)--(v) of \ref{mainthm2} are important for the rest of the paper, and we refer to them simply as (i)--(v).
Henceforth, we fix $\tau\ge 0$, and $\alpha, \delta\ge 1$, and $\beta\ge 2$, and some non-decreasing function
$\theta:\mathbb{N}\rightarrow \mathbb{N}$, for the remainder of the paper, and shall investigate the properties of
a graph satisfying (i)--(v).

Let $Y$ be a $(\zeta,\beta)$-core in $G$. A vertex $v\in V(G)\setminus Y$ is {\em dense}
to $Y$ if $v$ has at least $\alpha$ neighbours in each part of $Y$. We observe:
\begin{thm}\label{dense}
Let $G$ satisfy {\rm (i)--(v)},
and let $Y$ be a $(\zeta,\beta)$-core in $G$. Then there are at most $\alpha\tau2^{\beta\zeta}$ vertices in $G$ 
that are dense to $Y$.
\end{thm}
\Proof Let $A_1\ll A_{\beta}$ be the parts of $Y$,  and let $X_i\subseteq A_i$ with $|X_i|=\alpha$, for $1\le i\le \beta$. 
The set $N$ of vertices adjacent to all vertices in 
$X_1\cup\cdots\cup X_{\beta}$ has chromatic number at most $\tau$ by (ii), and includes 
no stable set of cardinality $\alpha$, since $G$ does not admit 
an $(\alpha,\beta+1)$-core by (v). Consequently $|N|\le \alpha\tau$. Since there are only at most $2^{\beta\zeta}$ choices for
$X_1\ll X_{\beta}$, and every vertex that is dense to $Y$ belongs to the set $N$ corresponding to some choice of
$X_1\ll X_{\beta}$, it follows that there are at most $\alpha\tau2^{\beta\zeta}$ vertices that are dense to $Y$.
This proves \ref{dense}.~\bbox

\section{Templates}

We will use an extension of the
template method of Kierstead-Penrice and Kierstead-Zhou, which was used in
different (and not easily compatible) ways in those papers.
Let $\eta\ge 1$ and $\zeta\ge \max(\eta,\alpha)$ be integers,
and let $G$ satisfy {\rm (i)--(v)}.
If $Y$ is a $(\zeta,\beta)$-core in $G$, we say that a vertex $v\in V(G)$ is {\em $\eta$-mixed} on $Y$ if 
\begin{itemize}
\item $v$ is not dense to $Y$; and
\item $v$ has at least $\eta$ neighbours in some part of $Y$.
\end{itemize}
Thus every vertex in $Y$ is $\eta$-mixed on $Y$.
A {\em $(\zeta,\eta)$-template} in $G$ is a pair $(Y,H)$, where $Y$ is a $(\zeta,\beta)$-core in $G$, and $H$ is a set of vertices
of $G$ with $Y\subseteq H$ such that every vertex in $H$ is $\eta$-mixed on $Y$. (Note that there may be vertices in 
$V(G)\setminus H$ that are $\eta$-mixed on $Y$.)

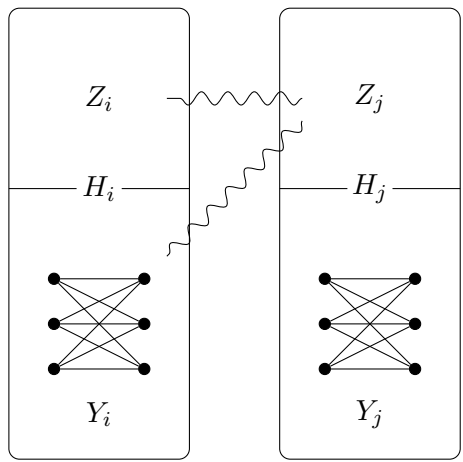
\begin{figure}[H]
\centering

\tikzset{snake it/.style={decorate, decoration=snake}}
\begin{tikzpicture}[scale=.6,auto=left]
\draw [rounded corners] (-5,-6) rectangle (-1,4);
\draw [rounded corners] (1,-6) rectangle (5,4);
\path[draw=black] (-5,0) -- (-3.5,0);
\path[draw=black] (-2.5,0) -- (-1,0);
\path[draw=black] (3.5,0) -- (5,0);
\path[draw=black] (1,0) -- (2.5,0);

\node at (-3,0) {$H_i$};
\node at (-3,-5) {$Y_i$};
\node at (3,-5) {$Y_j$};
\node at (3,0) {$H_j$};
\node at (-3,2) {$Z_i$};
\node at (3,2) {$Z_j$};

\tikzstyle{every node}=[inner sep=1.5pt, fill=black,circle,draw]
\node (a1) at (4,-4) {};
\node (a2) at (4,-3) {};
\node (a3) at (4,-2) {};
\node (b1) at (2,-4) {};
\node (b2) at (2,-3) {};
\node (b3) at (2,-2) {};
\node (c1) at (-2,-4) {};
\node (c2) at (-2,-3) {};
\node (c3) at (-2,-2) {};
\node (d1) at (-4,-4) {};
\node (d2) at (-4,-3) {};
\node (d3) at (-4,-2) {};

\foreach \from/\to in {a1/b1,a1/b2,a1/b3,a2/b1,a2/b2,a2/b3,a3/b1,a3/b2,a3/b3,c1/d1,c1/d2,c1/d3,c2/d1,c2/d2,c2/d3,c3/d1,c3/d2,c3/d3}
\draw [-] (\from) -- (\to);
\path[draw=black, snake it] (1.5,2) -- (-1.5,2);
\path[draw=black, snake it] (1.5,1.5) -- (-1.5,-1.5);

\end{tikzpicture}

\caption{Two terms of a $(\zeta,\eta)$-template sequence, with $j>i$. Wiggles indicate possible edges.} \label{fig:2}
\end{figure}

A {\em $(\zeta,\eta)$-template sequence} in $G$ is a sequence $(Y_i,H_i)\;(1\le i\le n)$ of $(\zeta,\eta)$-templates, 
such that
\begin{itemize}
\item for $1\le i<j\le n$, $H_i\cap H_j=\emptyset$; 
\item for $1\le i<j\le n$, there is no edge between $H_i$ and $Y_j$; and
\item for $1\le i<j\le n$, no vertex in $H_j$ is $\eta$-mixed on $Y_i$.
\end{itemize}
Later, we will denote $H_i\setminus Y_i$ by $Z_i$ (see figure \ref{fig:2}).

A {\em $(\zeta,\eta)$-template array} $\mathcal{T}$ in $G$ consists of a $(\zeta,\eta)$-template sequence $(Y_i,H_i)\;(1\le i\le n)$
together
with a set $U(\mathcal{T})\subseteq V(G)$, such that
for every vertex $v\in U(\mathcal{T})$, 
\begin{itemize}
\item $v\notin H_i$ for $1\le i\le n$; 
\item $v$ is not $\eta$-mixed on $Y_i$ for $1\le i\le n$; and
\item $v$ has a neighbour in $H_1\cup\cdots\cup H_n$.
\end{itemize}
We call $(Y_i,H_i)\;(1\le i\le n)$ the {\em sequence} of $\mathcal{T}$, and define 
$H(\mathcal{T}) = H_1\cup \cdots\cup H_n$ and 
$V(\mathcal{T})=H(\mathcal{T})\cup  U(\mathcal{T})$.

\begin{thm}\label{leftovers} 
Let $\eta\ge 1$ and $\zeta\ge \max(\eta,\alpha)$ be integers, and let $G$ satisfy {\rm (i)--(v)}.
Then there is a $(\zeta,\eta)$-template array $\mathcal{T}$ in $G$
such that $V(G)\setminus V(\mathcal{T})$ has
chromatic number at most $\theta(\zeta)$.
\end{thm}
\Proof
Let $(Y_i,H_i)\;(1\le i\le n)$ be a $(\zeta,\eta)$-template sequence with the property that for $1\le i\le n$, $H_i$ is the 
set of all vertices in $G$ that are $\eta$-mixed on $Y_i$, and subject to this, with $n$ maximum. 
Let $H=H_1\cup \cdots\cup H_n$, and let
$U$ be the set of vertices in $V(G)\setminus H$ with a neighbour in $H$. Let $\mathcal{T}$ be the $(\zeta,\eta)$-template array
consisting of $(Y_i,H_i)\;(1\le i\le n)$ together with $U(\mathcal{T})=U$. Let
$W=V(G)\setminus V(\mathcal{T})$, and suppose that there is a $(\zeta,\beta)$-core $Y_{n+1}\subseteq W$.
Let $H_{n+1}$ be the set of vertices in $G$ that are $\eta$-mixed on $Y_{n+1}$. Then $(Y_{n+1},H_{n+1})$ is a $(\zeta,\eta)$-template,
and no vertex in $H_{n+1}$ belongs to $H$, since no vertex in $H$ has a neighbour in $Y_{n+1}$, and every vertex in $H_{n+1}$
has a neighbour in $Y_{n+1}$. Consequently $(Y_i,H_i)\;(1\le i\le n+1)$ is a $(\zeta,\eta)$-template sequence, 
contrary to the maximality of $n$.
Thus there is no $(\zeta,\beta)$-core in $W$. From (iv),
$\chi(W)\le \theta(\zeta)$. This proves \ref{leftovers}.~\bbox

By setting $\phi(x)= x+\theta(\zeta)$ for $x\ge 0$, we deduce from \ref{leftovers} that:
\begin{thm}\label{getshadow}
Let $\eta\ge 1$ and $\zeta\ge \max(\eta,\alpha)$ be integers. 
Then there is a non-decreasing 
function $\phi:\mathbb{N}\rightarrow \mathbb{N}$ with the following property.
For all $c\ge 0$, if $G$ satisfies {\rm (i)--(v)} and $\chi(G)>\phi(c)$ then there is a $(\zeta,\eta)$-template array $\mathcal{T}$ 
in $G$ such that $\chi(V(\mathcal{T}))>c$.
\end{thm}

We will use the following elementary fact many times in the remainder of the paper, and we leave its proof to the reader:
\begin{thm}\label{digraph}
Let $D$ be a digraph with maximum outdegree at most $d$. Then the graph underlying $D$ has chromatic number at most $2d+1$,
and at most $d+1$ if $D$ is acyclic (that is, has no directed cycle).
\end{thm}

A $(\zeta,\eta)$-template array $\mathcal{T}$ with sequence $(Y_i,H_i)\;(1\le i\le n)$ is {\em partially 1-cleaned} if 
\begin{itemize}
\item for all distinct $i,j\in \{1\ll n\}$, no vertex of
$H_j$ is dense to $Y_i$; and
\item for all distinct $i,j\in \{1\ll n\}$ and all $v\in U(\mathcal{T})$, if $v$ is dense to $Y_i$ then $v$
has no neighbours in $H_j$;
\end{itemize}
and {\em 1-cleaned} if 
\begin{itemize}
\item for $1\le i\le n$, no vertex in $V(\mathcal{T})$ is dense to $Y_i$.
\end{itemize}

\begin{thm}\label{getpart1cleaned}
Let $\eta\ge 1$ and $\zeta\ge \max(\eta,\alpha)$ be integers.                               
Then there is a non-decreasing 
function $\phi:\mathbb{N}\rightarrow \mathbb{N}$ with the following property.
For all $c\ge 0$, if $G$ satisfies {\rm (i)--(v)} and $\chi(G)>\phi(c)$ then there is a partially 
1-cleaned $(\zeta,\eta)$-template array $\mathcal{T}$ 
in $G$ such that $\chi(V(\mathcal{T}))>c$.
\end{thm}
\Proof
Let $\psi$ satisfy \ref{getshadow} (with $\phi$ replaced by $\psi$). Let $t = \alpha\tau2^{\beta\zeta}$, and 
define $\phi(x) = \psi((2t+1)x)$ for all $x\in \mathbb{N}$; we claim that $\phi$ satisfies \ref{getpart1cleaned}. 
Let $c\ge 0$, and let $G$ satisfy (i)--(v), with
$\chi(G)>\phi(c)$.
By \ref{getshadow}, there is a $(\zeta,\eta)$-template array $\mathcal{T}$ 
in $G$ such that $\chi(V(\mathcal{T}))>(2t+1)c$. Let its sequence be $(Y_i,H_i)\;(1\le i\le n)$.
Choose a partition $U_i\;(1\le i\le n)$ of $U(\mathcal{T})$, such that 
\begin{itemize}
\item for $1\le i\le n$, every vertex of $U_i$ has a neighbour in $H_i$; and
\item for all distinct $i,j\in \{1\ll n\}$, if $v\in U_i$ is dense to $Y_i$ and $v$ has a neighbour in $H_j$
then $v$ is dense to $Y_j$.
\end{itemize}
(This can be arranged by assigning each vertex $v\in U(\mathcal{T})$ to some set $U_i$ where $v$ 
has a neighbour in $H_i$ and is not dense to $Y_i$ if possible, and otherwise assigning $v$
to some set $U_i$ where $v$ has a neighbour in $H_i$.)

Let $D$ be the digraph with vertex set $\{1\ll n\}$ in which for $1\le i,j\le n$ with $i\ne j$, if some vertex in $H_j\cup U_j$ is 
dense to $Y_i$
then $j$ is adjacent from $i$. 
By \ref{dense}, $D$ has maximum outdegree at most $t$, 
and so by \ref{digraph}
the graph underlying $D$ is 
$(2t+1)$-colourable. Consequently there is a partition $I_1\ll I_{2t+1}$ of $\{1\ll n\}$ such that for $1\le s\le {2t+1}$ 
and all distinct
$i,j\in I_s$, no vertex of $H_j\cup U_j$ is dense to $Y_i$. Hence the subsequence of $(Y_i,H_i)\;(1\le i\le n)$
consisting of the terms with $i\in I_s$, together with the set $\bigcup_{i\in I_s}U_i$, 
is a partially 1-cleaned $(\zeta,\eta)$-template array, $\mathcal{T}_s$ say. But every vertex
of $V(\mathcal{T})$ belongs to $V(\mathcal{T}_s)$ for some $s$; and so there
exists $s\in \{1\ll 2t+1\}$ such that 
$$\chi(V(\mathcal{T}_s))\ge \chi(V(\mathcal{T}))/(2t+1)>c.$$
This proves \ref{getpart1cleaned}.~\bbox

\begin{thm}\label{get1cleaned}
Let $\eta\ge 1$ and $\zeta\ge \max(\eta,\alpha)$ be integers.
Then there is a non-decreasing
function $\phi:\mathbb{N}\rightarrow \mathbb{N}$ with the following property.
For all $c\ge 0$, if $G$ satisfies {\rm (i)--(v)} and $\chi(G)>\phi(c)$ then there is a 
1-cleaned $(\zeta,\eta)$-template array $\mathcal{T}$
in $G$ such that $\chi(V(\mathcal{T}))>c$.
\end{thm}
\Proof
Let $\psi$ satisfy \ref{getpart1cleaned} (with $\phi$ replaced by $\psi$). 
Let $t = \alpha\tau2^{\beta\zeta}$, and define 
$\phi(x) = \psi(x+t\tau)$ for all $x\in \mathbb{N}$; we claim that $\phi$ satisfies \ref{get1cleaned}.
Let $c\ge 0$, and let $G$ satisfy (i)--(v), with
$\chi(G)>\phi(c)$. By \ref{getpart1cleaned}, there is a partially
1-cleaned $(\zeta,\eta)$-template array $\mathcal{T}$
in $G$ such that $\chi(V(\mathcal{T}))>c+t\tau$. Let its sequence be $(Y_i,H_i)\;(1\le i\le n)$;
and for $1\le i\le n$ let $X_i$ be the set of all vertices in $V(\mathcal{T})$ that are
dense to $Y_i$. Hence $X_i\subseteq U(\mathcal{T})$, and for all distinct $i,j\in \{1\ll n\}$,
no vertex of $X_i$ has a neighbour in $H_j$. Let $X=X_1\cup\cdots\cup X_n$.
By \ref{dense} $|X_i|\le t$ for each $i$.
Hence we may partition $X$ into $t$ sets $W_1\ll W_t$ such that $|X_i\cap W_j|\le 1$
for all $i,j$ with $1\le i\le n$ and $1\le j\le t$. For $v\in W_j$, let $v\in X_i$; then $v$ has a neighbour in $Y_i$,
and has no neighbours in $Y_{i'}$ for $i'\ne i$; and so each set $W_j$ is matching-covered in $G$.
By (iii), $\chi(W_j)\le \tau$ for each $j$, and so $\chi(W_1\cup\cdots\cup W_t)\le t\tau$.
This proves that $\chi(X)\le t\tau$. Let $\mathcal{T}'$ be the $(\zeta,\eta)$-template array in $G$ with the same sequence as
$\mathcal{T}$ and with $U(\mathcal{T}') = U(\mathcal{T}\setminus X)$; then $\mathcal{T}'$ is 1-cleaned, and 
$$\chi(V(\mathcal{T'}))\ge \chi(V(\mathcal{T}))-t\tau>c.$$
This proves \ref{get1cleaned}.~\bbox

\section{Edges between templates}

\begin{thm}\label{manyY} 
Let $\eta\ge 1$ and $\zeta\ge \max(\eta+\delta,\alpha)$ be integers, and let $G$ satisfy {\rm (i)--(v)}.
Let $\mathcal{T}$ be a 1-cleaned $(\zeta,\eta)$-template array in $G$, with sequence  $(Y_i,H_i)\;(1\le i\le n)$.
For each $v\in V(\mathcal{T})$, there are at most $2\delta$ values of $i\in \{1\ll n\}$ such that $v$ has a neighbour in $Y_i$.
\end{thm}
\Proof
Suppose there exists $I\subseteq \{1\ll n\}$ with $|I|=2\delta+1$ such that $v$ has a neighbour in $Y_i$ for each $i\in I$,
and let $i_0$ be the maximum element of $I$. Let $I'=I\setminus \{i_0\}$.
It follows that $v\notin H_i$ for all $i\in I'$. Consequently for $i\in I'$, $v$ is not $\eta$-mixed on $Y_i$,
and not dense to $Y_i$ (since the template array is 1-cleaned and $v\in V(\mathcal{T})$); and it follows that 
$v$ has at most $\eta-1$ neighbours in each part of $Y_i$. Let $i\in I'$, and let the parts of $Y_i$ be $A_1\ll A_{\beta}$.
Since $v$ has a neighbour in $Y_i$, we may assume that $v$ has a neighbour in $A_1$. Since $v$ has at most $\eta-1$
neighbours in $A_2$, and $|A_2|=\zeta$, it follows that $v$ has at least $\delta$ non-neighbours in $A_2$, and so there is a 
$(1,\delta)$-broom with handle $v$ in $G[Y_i\cup \{v\}]$. 
But also, since $v$ has at least $\delta$ non-neighbours in $A_1$, 
there is a $(2,\delta)$-broom with handle $v$ in $G[Y_i\cup \{v\}]$. By selecting the $(1,\delta)$-broom with handle $v$ 
in $G[Y_i\cup \{v\}]$ for $\delta$ values of $i\in I'$, and selecting the $(2,\delta)$-broom 
for the remaining $\delta$ values of $i\in I'$, and taking the union of all these brooms, we find that $G$
contains $T(\delta)$, contrary to (i). This proves \ref{manyY}.~\bbox

For the remainder of the paper, let us define $\gamma= (2\delta\tau+1)(2\delta+1)$.
\begin{thm}\label{manyH} 
Let $\eta\ge \delta$ and $\zeta\ge \max(\eta,\alpha)+\delta$ be integers. 
Let $G$ satisfy {\rm (i)--(v)},
and let $\mathcal{T}$ be a 1-cleaned 
$(\zeta,\eta)$-template array in $G$, with sequence $(Y_i,H_i)\;(1\le i\le n)$.
For each $v\in V(\mathcal{T})$, there are fewer than $\gamma$ values of $i\in \{1\ll n\}$ such that $v$ has a neighbour in $H_i$.
\end{thm}
\Proof 
Suppose then that $G$, $\mathcal{T}$, and $(Y_i,H_i)\;(1\le i\le n)$ are as in the theorem, and there are at least $\gamma$ values of 
$i\in \{1\ll n\}$ such that $v$ has a neighbour in $H_i$.
By \ref{manyY}, there are at most $2\delta+1$ values of $i$ such that $v$ has a neighbour in $Y_i$.
Consequently there exists 
$I_1\subseteq \{1\ll n\}$ with $|I_1|= 2\delta(2\delta+1)\tau$ 
such that for each $i\in I_1$,
$v$ has a neighbour in $H_i$ and $v$ has no neighbour in $Y_i$.
For $i\in I_1$, let $u_i \in H_i\setminus Y_i$ be adjacent to $v$.

Let $D$ be the digraph with vertex set $I_1$ in which for distinct $i,j\in I_1$, $i$ is adjacent from $j$ in $D$ if $u_j$
has a neighbour in $Y_i$ (and consequently $i<j$). From \ref{manyY}, $D$ has maximum outdegree at most $2\delta$ 
and is acyclic, and so by \ref{digraph} the graph underlying
$D$ has chromatic number at most $2\delta+1$. Hence there exists $I_2\subseteq I_1$ with 
$|I_2|= |I_1|/(2\delta+1)=2\delta\tau$
such that for all distinct $i,j\in I_2$, $u_j$ has no neighbour in $Y_i$.
Since $\chi(\{u_i:i\in I_2\}\le \tau$ by (ii),
there exists $I_3\subseteq I_2$ with $|I_3|= 2\delta$
such that for all $i<j$ with $i,j\in I_3$, $u_i$ and $u_j$ are nonadjacent. For each $i\in I_3$, since $u_i$ is $\eta$-mixed on $Y_i$,
and $\eta\ge \delta$, and $v$ has no neighbour in $Y_i$, it follows that
there is a $(1,\delta)$-broom in $G[\{v,u_i\}\cup Y_i]$ with handle $v$. Since $u_i$ is $\eta$-mixed on $Y_i$,
and has fewer than $\alpha$ neighbours in some part of $Y_i$, and $\zeta\ge \delta+\alpha$, it follows that
there are two distinct parts $A_1,A_2$ of $Y_i$ such that $u_i$ has a neighbour in $A_1$
and has at least $\delta$ non-neighbours in $A_2$.
Consequently there is a 
$(2,\delta)$-broom in $G[\{v,u_i\}\cup Y_i]$ with handle $v$. But then, choosing the $(1,\delta)$-broom for $\delta$
values of $i\in I_3$ and choosing the $(2,\delta)$-broom for the other $\delta$
values of $i\in I_3$, and taking their union, we find that $G$ contains $T(\delta)$, contrary to (i). 
This proves \ref{manyH}.~\bbox

\begin{thm}\label{manystrongH}
Let $\eta\ge \delta$ and $\zeta\ge \max(\eta,\alpha)+\delta$ be integers. Then there exists
$s$ with the following property.
Let $G$ satisfy {\rm (i)--(v)},
and let $\mathcal{T}$ be a 1-cleaned
$(\zeta,\eta)$-template array in $G$, with sequence $(Y_i,H_i)\;(1\le i\le n)$.
For $1\le j\le n$, there are at most $s$ values of $i\in \{1\ll n\}$ such that some vertex in $H_j$ has
at least $\delta$ neighbours in $H_i$.
\end{thm}
\Proof
Let $s_3=2\delta\tau$, let $s_2= (2(\delta+1)\gamma+1)s_3$, let $s_1= s_2+\gamma$, and 
let $s=\zeta\beta s_1$.
Now let $\eta,\zeta,G$, $\mathcal{T}$ and $(Y_i,H_i)\;(1\le i\le n)$
be as in the theorem, and suppose that for some $j\in \{1\ll n\}$ and some subset $I\subseteq \{1\ll n\}$ with $|I|>s$,
and every $i\in I$, there exists $u_i\in H_j$ with at least $\delta$ neighbours in $H_i$. Since each $u_i$ has a neighbour in $Y_j$,
and $|Y_j| = \zeta\beta$, there exists $I_1\subseteq I$ with $|I_1|= s_1$ and a vertex $y\in Y_j$
adjacent to every $u_i\;(i\in I_1)$. Since by \ref{manyH}, $y$ has neighbours in $H_i$ for at most $\gamma$ values of $i$,
there exists $I_2\subseteq I_1$ with $|I_2|= |I_1|-\gamma=s_2$ such that $y$ has no neighbours in $H_i$ for $i\in I_2$ 
(and in particular $j\notin I_2$). For each $i\in I_2$, choose a set $W_i\subseteq H_i$ with $|W_i|=\delta$ such that every vertex in $W_i$
is adjacent to $u_i$. 

Let $D$ be the digraph with vertex set $I_2$ in which for distinct $i,i'\in I_1$,
$i'$ is adjacent from $i$ if some vertex in $u_{i}\cup W_{i}$ has a neighbour in $H_{i'}$. 
(Thus if $u_i = u_{i'}i$ then we have $i\rightarrow i'$ and $i'\rightarrow i$.) By \ref{manyH}, $D$
has maximum outdegree at most $(\delta+1)\gamma$, and so by \ref{digraph} the graph underlying $D$ has chromatic number at most $2(\delta+1)\gamma+1$.
Hence there exists $I_3\subseteq I_2$ with $|I_3|= |I_2|/(2(\delta+1)\gamma+1)= s_3$ such that for all distinct $i,i'\in S_3$,
no vertex in $u_{i}\cup W_{i}$ has a neighbour in $H_{i'}$ (and in particular the vertices $u_i\;(i\in I_3)$ are all
distinct). By (ii) the set $\{u_i:i\in I_3\}$ has chromatic number at most $\tau$, so
there exists $I_4\subseteq I_3$ with $|I_4|= 2\delta$
such that the vertices $u_i\;(i\in I_4)$ are pairwise nonadjacent. 

For each $i\in I_4$, there is a $(1,\delta)$-broom
in $G[\{y,u_i\}\cup W_i]$ with handle $y$. Moreover, choose $w_i\in W_i$; then since $w_i$ is $\eta$-mixed on $Y_i$,
there is a $(2,\delta)$-broom in  $G[\{y,u_i, w_i\}\cup Y_i]$ with handle $y$. By choosing the $(1,\delta)$-broom
for $\delta$ values of $i$, and the $(2,\delta)$-broom for the other $\delta$ values of $i$, and taking their union,
we find that $G$ contains $T(\delta)$, contrary to (i). This proves \ref{manystrongH}.~\bbox

For $d\ge 0$, let us say 
a 1-cleaned $(\zeta,\eta)$-template array $\mathcal{T}$ with sequence $(Y_i,H_i)\;(1\le i\le n)$ is {\em partially $(2,d)$-cleaned} if
\begin{itemize}
\item for all $i\in \{1\ll n\}$, every vertex of $H_i$ has at most $d$ neighbours in $H(\mathcal{T})\setminus H_i$
\end{itemize}
and {\em 2-cleaned} if 
\begin{itemize}
\item for all distinct $i,j\in \{1\ll n\}$, no vertex of $H_i$ has a neighbour in $H_j$, and
\item for $1\le i\le n$, $H_i\setminus Y_i$ is stable.
\end{itemize}

\begin{thm}\label{getpart2cleaned}
Let $\eta\ge \delta$ and $\zeta\ge \max(\eta,\alpha)+\delta$ be integers.
Then there exist $d\ge 0$ and a non-decreasing
function $\phi:\mathbb{N}\rightarrow \mathbb{N}$, with the following property.
For all $c\ge 0$, if $G$ satisfies {\rm (i)--(v)} and $\chi(G)>\phi(c)$ then there is a
partially $(2,d)$-cleaned $(\zeta,\eta)$-template array $\mathcal{T}$
in $G$ such that $\chi(V(\mathcal{T}))>c$.
\end{thm}
\Proof
Let $s$ be as in \ref{manystrongH}, and let $d=\gamma(\delta-1)$.
Let $\psi$ satisfy \ref{get1cleaned} (with $\phi$ replaced by $\psi$), and define
$\phi(x) = \psi((2s+1)x)$ for $x\ge 0$. Now let $c\ge 0$, and let $G$ satisfy {\rm (i)--(v)}, with $\chi(G)>\phi(c)$.
By \ref{get1cleaned} there is 
a 1-cleaned $(\zeta,\eta)$-template array $\mathcal{T}$
in $G$ such that $\chi(V(\mathcal{T}))>(2s+1)c$.
Let its sequence be $(Y_i,H_i)\;(1\le i\le n)$.
Let $D$ be the digraph with vertex set $\{1\ll n\}$ in which for distinct $i,j$ with $1\le i,j\le n$,
$i$ is adjacent from $j$ if some vertex of $H_j$ has at least $\delta$ neighbours in $H_i$. By \ref{manystrongH},
$D$ has maximum outdegree at most $s$, and so by \ref{digraph} the graph underlying $D$ has chromatic number at most $2s+1$.
Consequently there is a partition $I_1\ll I_{2s+1}$ of $\{1\ll n\}$ such that for $1\le r\le 2s+1$, if
$i,j\in I_r$ are distinct then each vertex of $H_j$ has at most $\delta-1$ neighbours in $H_i$. By \ref{manyH} it follows that
for each $j\in I_r$, each vertex of $H_j$ has at most $\gamma(\delta-1)=d$ neighbours in $\bigcup_{i\in I_r} H_i\setminus H_j$.
For each $r\in \{1\ll 2s+1\}$, let $\mathcal{T}_r$ be the $(\zeta,\eta)$-template array with sequence the subsequence
of $(Y_i,H_i)\;(1\le i\le n)$ consisting of the terms with $i\in I_r$, and with $U(\mathcal{T}_r)$ the set of vertices
in $U(\mathcal{T})$ with a neighbour in $\bigcup_{i\in I_r} H_i$. Thus each $\mathcal{T}_r$ is partially 
$(2,d)$-cleaned; and since every vertex of $V(\mathcal{T})$ belongs to $V(\mathcal{T}_r)$ for some $r$, there exists
$r\in \{1\ll 2s+1\}$ such that $\chi(V(\mathcal{T}_r))\ge \chi(V(\mathcal{T}_r))/(2s+1)>c$. This proves \ref{getpart2cleaned}.~\bbox

\begin{thm}\label{get2cleaned}
Let $\eta\ge \delta$ and $\zeta\ge \max(\eta,\alpha)+\delta$ be integers.
Then there is a non-decreasing
function $\phi:\mathbb{N}\rightarrow \mathbb{N}$ with the following property.
For all $c\ge 0$, if $G$ satisfies {\rm (i)--(v)} and $\chi(G)>\phi(c)$ then there is a
$2$-cleaned $(\zeta,\eta)$-template array $\mathcal{T}$
in $G$ such that $\chi(U(\mathcal{T}))>c$.
\end{thm}
\Proof
Let $d, \psi$ satisfy \ref{getpart2cleaned} (with $\phi$ replaced by $\psi$).
Let $t=(d+1)\beta\zeta\tau$, and define
define $\phi(x) = \psi((x+\beta+1)t)$ for $x\ge 0$; we claim this satisfies the theorem. For let $c\ge 0$, and let $G$ satisfy {\rm (i)--(v)}
with $\chi(G)>\phi(c)$. By \ref{getpart2cleaned} 
there is a
partially $(2,d)$-cleaned $(\zeta,\eta)$-template array $\mathcal{T}$
in $G$ such that $\chi(V(\mathcal{T}))>(c+\beta+1)t$.
Let its sequence be $(Y_i,H_i)\;(1\le i\le n)$. For $1\le i\le n$, every vertex of $H_i$ has a neighbour
in $Y_i$, and since $|Y_i|=\beta\zeta$, and by (ii) the set of vertices in $G$ adjacent to any given vertex of $Y_i$
has chromatic number at most $\tau$, it follows that $\chi(H_i)\le \beta\zeta\tau$. 
Let $J_1$ be the subgraph of $G$ with vertex set $H(\mathcal{T})$ and edge set all edges of $G$
with an end in $H_i$ and an end in $H_j$ for distinct $i,j$; and let $J_2$ be the subgraph of $G$ with vertex set
$H(\mathcal{T})$ and edge set all edges of $G$ with both ends in $H_i$ for some $i$. We have just seen that $J_2$
has chromatic number at most $\beta\zeta\tau$; and since $\mathcal{T}$ is partially $(2,d)$-cleaned,
$J_1$ has maximum degree at most $d$ and so is $(d+1)$-colourable. Hence $G_1\cup G_2$ has chromatic number at
most $(d+1)\beta\zeta\tau=t$. Consequently there is a partition $W_1\ll W_t$ of $H_1\cup\cdots\cup H_n$ into $t$ stable sets.
For $1\le j\le t$, let $\mathcal{T}_j$ be the $(\zeta,\eta)$-template
with sequence $(Y_i,(H_i\cap W_j)\cup Y_i)\;(1\le i\le n)$, where $U(\mathcal{T}_j)$ is the set of vertices in
$U(\mathcal{T})$ with a neighbour in $W_j$. Then each $\mathcal{T}_j$ is 2-cleaned, and since every vertex of $V(\mathcal{T})$
belongs to $V(\mathcal{T}_j)$ for some $j$, there exists $j\in \{1\ll t\}$ such that 
$$\chi(V(\mathcal{T}_j))\ge \chi(V(\mathcal{T}))/t> c+\beta+1.$$ 
But $\chi(H(\mathcal{T}_j))\le \beta+1$, since 
$\mathcal{T}_j$ is 2-cleaned; and so $\chi(U(\mathcal{T}_j)) >c$. This proves \ref{get2cleaned}.~\bbox

\section{Shadowing, and growing daisies}

Let $\mathcal{T}$ be a
$(\zeta,\eta)$-template array in $G$, with sequence $(Y_i,H_i)\;(1\le i\le n)$.
A {\em shadowing} of $\mathcal{T}$ is a sequence $B_1\ll B_n$ of pairwise disjoint subsets of $U(\mathcal{T})$,
with union $U(\mathcal{T})$, such that for $1\le i\le n$, every vertex in $B_i$ has a neighbour in $H_i$.
Every template array has a shadowing, and in general it has many. Let us say a shadowing $B_1\ll B_n$ has {\em degree at most $s$}
if for every vertex $v\in V(\mathcal{T})$, there are at most $s$ values of $i\in \{1\ll n\}$ such that $v$ has a neighbour
in $B_i$. If $X\subseteq U(\mathcal{T})$, we say the shadowing has {\em degree at most $s$ relative to $X$} if 
for every vertex $v\in V(\mathcal{T})$, there are at most $s$ values of $i\in \{1\ll n\}$ such that $v$ has a neighbour
in $B_i\cap X$.

Let $B_1\ll B_n$ be a shadowing of $\mathcal{T}$. A {\em daisy}  (with respect to $\mathcal{T}$ and the given shadowing)
is an induced subgraph $D$ of $G$ 
isomorphic to the $\delta+2$-vertex star $K_{1,\delta+1}$, such that 
\begin{itemize}
\item exactly one vertex $u$ of $D$ belongs to $H(\mathcal{T})$; let $u\in H_i$ say;
\item $u$ has degree one in $D$, and the neighbour $v$ of $u$ in $D$ belongs to $U(\mathcal{T})$; 
\item there exists $j\ne i$ with $1\le j\le n$ such that $V(D)\setminus \{u,v\}\subseteq B_j$.
\end{itemize}
We call $u$ the {\em root}, $v$ the {\em eye}, and the vertices in $V(D)\setminus \{u,v\}$ the {\em petals}
of the daisy.
We need the following, proved in \cite{cuttingedge}, but we repeat the proof because it is short:
\begin{thm}\label{creatures}
Let $d\ge 0$ be an integer, let $G$ be a graph with chromatic number more than $d$, and let $X\subseteq V(G)$ be stable,
such that $\chi(G\setminus X)<\chi(G)$. Then some vertex in $X$ has at least $d$ neighbours in $V(G)\setminus X$.
\end{thm}
\Proof
Let $\chi(G)=k+1$, and so $k\ge d$. Let $\phi:V(G)\setminus X\rightarrow \{1\ll k\}$ be a $k$-colouring of $G\setminus X$.
For each $x\in X$,
if $x$ has at most $d-1$ neighbours in $V(G)\setminus X$ then we may choose $\phi(x)\in \{1\ll k\}$, different from $\phi(v)$
for each neighbour $v\in V(G)\setminus X$ of $x$; and this extends $\phi$ to a $k$-colouring of $G$, which is impossible. Thus
for some $x\in X$, $x$ has at least $d$ neighbours in $V(G)\setminus X$.
This proves \ref{creatures}.~\bbox

We deduce:
\begin{thm}\label{firstdaisy}
Let $\mathcal{T}$ be a
$(\zeta,\eta)$-template array in $G$, with sequence $(Y_i,H_i)\;(1\le i\le n)$.
Let $X\subseteq U(\mathcal{T})$ with $\chi(X)>s\beta\delta\zeta\tau^2$.
Let $B_1\ll B_n$ be a shadowing of degree at most $s$ relative to $X$.
Then there is a daisy in $G[H(\mathcal{T})\cup X]$.
\end{thm}
\Proof
For $1\le i\le n$, let $W_i$ be the set of vertices in $X$ that have a neighbour in $H_i$ and have none in 
$H_1\cup\cdots\cup H_{i-1}$.
Since every vertex in $W_i$ has distance at most two from some vertex in $Y_i$, and $|Y_i|=\beta\zeta$, it follows that
$\chi(W_i)\le \beta\zeta\tau$ for each $i$. Consequently there is a partition $X_1\ll X_{k\beta\tau}$ of $X$
such that $W_i\cap X_j$ is stable for $1\le i\le n$ and $1\le j\le \beta\zeta\tau$. We may assume that 
$\chi(X_1)\ge \chi(X)/(\beta\zeta\tau)> s\delta\tau$.
Choose $i$ minimum such that $\chi(X_1\setminus (W_1\cup\cdots W_{i}))<\chi(X_1)$, and let 
$H= G[X_1\setminus (W_1\cup\cdots W_{i-1})]$. Then since 
$\chi(H\setminus (X_1\cap W_{i}))< \chi(H)$, and $X_1\cap W_i$ is stable, and $\chi(H)> s\delta\tau$, 
\ref{creatures} implies that
some vertex $v\in W_i$ has a set $P$ of at least $s\delta\tau$ neighbours in $W_{i+1}\cup \cdots\cup W_n$. Choose $u\in H_i$
adjacent to $v$. Since the vertices in $P$
have no neighbours in $H_i$, they are nonadjacent to $u$. By hypothesis there are at most $s$ 
values of $j\in \{1\ll n\}$ such that $v$ has a neighbour in $B_j\cap X$, and so there exists $j\in \{1\ll n\}$
such that $|P\cap B_j|\ge \delta\tau$. Since $\chi(P)\le \tau$ (because the vertices in $P$ are all adjacent to $v$), 
there is a stable subset $P'\subseteq P\cap B_j$ with $|P'|\ge \delta$. Now $j\ne i$
since no vertices in $P$ have a neighbour in $H_i$; and so $G[\{u,v\}\cup P']$ is a daisy. This proves \ref{firstdaisy}.~\bbox

\bigskip

Let $\mathcal{T}$ be a
$(\zeta,\eta)$-template array in $G$, with sequence $(Y_i,H_i)\;(1\le i\le n)$; and 
let $B_1\ll B_n$ be a shadowing. A {\em bunch} of daisies is a set $\{D_j:j\in J\}$ of daisies 
where $J\subseteq \{1\ll n\}$ and 
for each $j\in J$, $D_j$ has root $u_j$, eye $v_j$ and set of petals $P_j$, with the following properties:
\begin{itemize}
\item $P_j\subseteq B_j$ for each $j\in J$;
\item there exists $i\in \{1\ll n\}\setminus J$ such that $u_j\in H_i$ for each $j\in J$;
\item for all distinct $j,j'\in J$, $P_j\cup \{v_j\}$ is disjoint from $P_{j'}\cup \{v_{j'}\}$, 
and there is no edge joining these two sets; and
\item for all distinct $j,j'\in J$, $u_j$ has no neighbour in $P_{j'}$.
\end{itemize}
(Thus, the roots may not all be distinct, and the root of one daisy may be adjacent to the eye of another.)
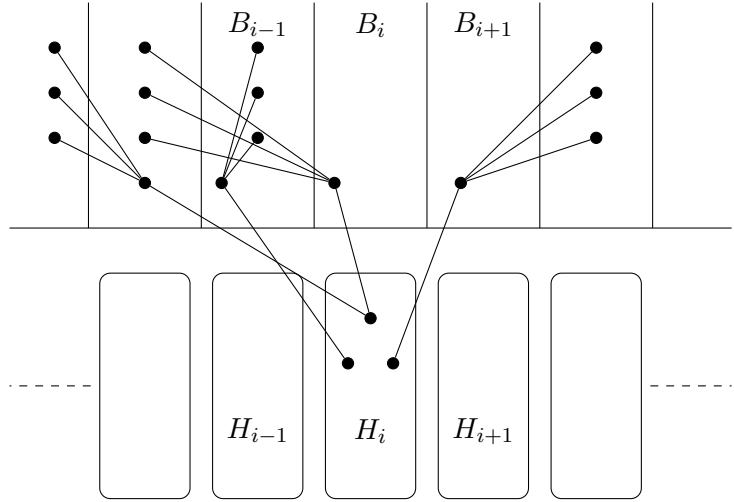
\begin{figure}[H]
\centering

\begin{tikzpicture}[scale=.6,auto=left]
\draw [rounded corners] (1,-6) rectangle (-1,-1);
\draw [rounded corners] (3.5,-6) rectangle (1.5,-1);
\draw [rounded corners] (-1.5,-6) rectangle (-3.5,-1);
\draw [rounded corners] (6,-6) rectangle (4,-1);
\draw [rounded corners] (-6,-6) rectangle (-4,-1);

\path[draw=black] (8,0) -- (-8,0);
\path[dashed, draw=black] (6.2,-3.5) -- (8,-3.5);
\path[dashed, draw=black] (-6.2,-3.5) -- (-8,-3.5);

\node at (0,-4.5) {$H_i$};
\node at (2.5,-4.5) {$H_{i+1}$};
\node at (-2.5,-4.5) {$H_{i-1}$};
\node at (0,4.5) {$B_i$};
\node at (2.5,4.5) {$B_{i+1}$};
\node at (-2.5,4.5) {$B_{i-1}$};

\path[draw=black] (1.25,0) -- (1.25,5);
\path[draw=black] (-1.25,0) -- (-1.25,5);
\path[draw=black] (3.75,0) -- (3.75,5);
\path[draw=black] (-3.75,0) -- (-3.75,5);
\path[draw=black] (6.25,0) -- (6.25,5);
\path[draw=black] (-6.25,0) -- (-6.25,5);

\tikzstyle{every node}=[inner sep=1.5pt, fill=black,circle,draw]

\node (u1) at (0,-2) {};
\node (u2) at (-0.5,-3) {};
\node (u3) at (0.5,-3) {};
\node (v1) at (-0.8,1) {};
\node (v2) at (-3.3,1) {};
\node (v3) at (2,1) {};
\node (v4) at (-5,1) {};
\node (a1) at (-5,2) {};
\node (a2) at (-2.5,2) {};
\node (a3) at (5,2) {};
\node (a4) at (-7,2) {};
\node (b1) at (-5,3) {};
\node (b2) at (-2.5,3) {};
\node (b3) at (5,3) {};
\node (b4) at (-7,3) {};
\node (c1) at (-5,4) {};
\node (c2) at (-2.5,4) {};
\node (c3) at (5,4) {};
\node (c4) at (-7,4) {};

\foreach \from/\to in {u1/v1, u2/v2,u3/v3,u1/v4,v1/a1,v1/b1,v1/c1,v2/a2,v2/b2,v2/c2,v3/a3,v3/b3,v3/c3,v4/a4,v4/b4,v4/c4}
\draw [-] (\from) -- (\to);

\end{tikzpicture}

\caption{A bunch of daisies.} \label{fig:3}
\end{figure}

We deduce:
\begin{thm}\label{getdaisies}
Let $\mathcal{T}$ be a
$(\zeta,\eta)$-template array in $G$, with sequence $(Y_i,H_i)\;(1\le i\le n)$.
Let $t\ge 0$, and let $X\subseteq U(\mathcal{T})$ with 
$$\chi(X)> 2st\zeta^2\delta(2(\delta+2)s+3)\beta^2\tau^3.$$
Let $B_1\ll B_n$ be a shadowing of degree at most $s$ relative to $X$.
Then there is a bunch $\{D_1\ll D_t\}$ of daisies, such that $V(D_j)\subseteq H(\mathcal{T})\cup X$ for $1\le j\le t$.
\end{thm}
\Proof
Let $m_1=t\beta\zeta\tau$, and $m=t\beta\zeta\tau(2(\delta+2)s+3)$.
Let $D$ be the digraph with vertex set $\{1\ll n\}$ in which for distinct $i,j\in \{1\ll n\}$, $j$ is adjacent from $i$
if there is a daisy with root in $H_i$, eye in $X$ and set of petals in $X\cap B_j$. 
\\
\\
(1) {\em There exists $i\in \{1\ll n\}$ with outdegree in $D$ at least $m$.}
\\
\\
Suppose not; then by \ref{digraph} the graph underlying $D$ is $2m$-colourable. Consequently there exists $I\subseteq \{1\ll n\}$
with 
$$\chi(X')\ge \chi(\bigcup_{1\le i\le n}B_i\cap X)/(2m)> s\delta\beta\zeta\tau^2$$ 
where $X'= \bigcup_{i\in I}B_i\cap X$,
such that for all distinct $i,j\in I$, there is no daisy with root in $H_i$, eye in $X$ and set of petals in $X\cap B_j$.
In particular, applying \ref{firstdaisy} to the $(\zeta,\eta)$-template
array $\mathcal{T}'$ with sequence $(Y_i,H_i)(i\in I)$ and $U(\mathcal{T}')=X'$, it follows that
$\chi(X')\le s\beta\delta\zeta\tau^2$, a contradiction. This proves (1).

\bigskip

From (1), there exist $i\in \{1\ll n\}$ and $J\subseteq \{1\ll n\}\setminus \{i\}$, with $|J|=m$, such that
for each $j\in J$ there is a daisy $D_j$ with root $u_j\in H_i$, eye $v_j\in X$ and set of petals $P_j\subseteq B_j\cap X$.
Now let $D'$ be the digraph with vertex set $J$ in which for all distinct $j,j'\in J$, $j$ is adjacent from $j'$ if some vertex in 
$D_j$ belongs to or has a neighbour in $P_{j'}$. Then $D'$ has maximum outdegree at most $(\delta+2)s+1$, 
and so  by \ref{digraph} the graph underlying $D'$ is
$2(\delta+2)s+3$-colourable. Hence there exists $J_1\subseteq J$ with $|J_1|= m/(2(\delta+2)s+3)=m_1$, such that 
for all distinct $j,j'$, no vertex in $D_j$ has a neighbour in $P_{j'}$. In particular, the vertices 
$v_j(j\in J_1)$ are all distinct. Since $\{v_j:j\in J_1\}$ has chromatic number at most $\beta\zeta\tau$, there exists 
$J_2\subseteq J_1$ with $|J_2|=m_1/(\beta\zeta\tau) = t$ such that the vertices $v_j(j\in j_2)$ are pairwise nonadjacent.
But then $\{D_j:j\in J_2\}$ is a bunch of daisies of cardinality $t$. This proves \ref{getdaisies}.~\bbox

\section{Privatization}

Let $A,B$ be disjoint subsets of $V(G)$; we say $A$ {\em covers} $B$ if every vertex in $B$ has a
neighbour in $A$. 
We claim:
\begin{thm}\label{privatecover}
Let $A,B\subseteq V(G)$ be disjoint, and let $A$ cover $B$. Let $d\ge 0$ be an integer.
Then there exist $A'\subseteq A$ and $B'\subseteq B$ such that
\begin{itemize}
\item $A'$ covers $B\setminus B'$;
\item $B'$ is the union of $d$ matching-covered sets;
\item every vertex in $B'$ has at most one neighbour in $A'$; and
\item every vertex in $A'$ has exactly $d$ neighbours in $B'$.
\end{itemize}
\end{thm}
\Proof
We proceed by induction on $d$. The result is trivial for $d=0$, because we can set $A'=A$ and $B'=\emptyset$;
so we assume that $d>0$ and the result
holds for $d-1$. Hence there exist $A'\subseteq A$ and $B''\subseteq B$ such that
\begin{itemize}
\item $A'$ covers $B\setminus B''$;
\item $B''$ is the union of $d-1$ matching-covered sets;
\item every vertex in $B''$ has at most one neighbour in $A'$; and
\item every vertex in $A'$ has exactly $d-1$ neighbours in $B''$.
\end{itemize}
Choose $A'$ minimal with this property. Consequently for each $u\in A'$ there exists $v_u\in B\setminus B''$ such that $u$
is the unique neighbour of $v_u$ in $A'$. Let $X=\{v_u:u\in A'\}$ and let $B'=B''\cup X$.
Then $X$ is matching-covered, and every vertex in $A'$ has a unique neighbour in $X$ and exactly $d-1$ in $B''$,
and so exactly $d$ in $B'$. Consequently $B'$ satisfies the theorem. This proves \ref{privatecover}.~\bbox

Let $\mathcal{T}$ be a 2-cleaned
$(\zeta,\eta)$-template array in $G$, with sequence $(Y_i,H_i)\;(1\le i\le n)$. We recall that
$H(\mathcal{T})$ denotes $\bigcup_{1\le i\le n}H_i$; and let $Y(\mathcal{T})$ denote $\bigcup_{1\le i\le n}Y_i$
and $Z(\mathcal{T})=H(\mathcal{T})\setminus Y(\mathcal{T})$.
A {\em privatization} for $\mathcal{T}$
is a subset $\Pi\subseteq U(\mathcal{T})$ such that 
\begin{itemize}
\item $\Pi$ is the union of $\delta\tau$ matching-covered sets;
\item every vertex in $\Pi$ has exactly one neighbour in $Z(\mathcal{T})$ and none in $Y(\mathcal{T})$; and
\item every vertex in $Z(\mathcal{T})$ has exactly $\delta\tau$ neighbours in $\Pi$.
\end{itemize}
We deduce:
\begin{thm}\label{privatize}
Let $G$ satisfy {\rm (i)--(v)}, and 
let $\mathcal{T}$ be a 2-cleaned
$(\zeta,\eta)$-template array in $G$, with sequence $(Y_i,H_i)\;(1\le i\le n)$. Then there
is a 2-cleaned $(\zeta,\eta)$-template array $\mathcal{T}'$, with sequence $(Y_i,H_i')\;(1\le i\le n)$
and a privatization $\Pi$ for $\mathcal{T}'$,
such that 
\begin{itemize}
\item $H_i'\subseteq H_i$ for $1\le i\le n$; 
\item $U(\mathcal{T}')\subseteq U(\mathcal{T})$; and
\item $\chi(U(\mathcal{T}')\setminus \Pi)\ge \chi(U(\mathcal{T}))-\delta\tau^2$.
\end{itemize}
\end{thm}
\Proof We will obtain the desired $\mathcal{T}'$ by removing some elements of each $H_i$ and also removing some elements of 
$U(\mathcal{T})$. We cannot remove from $H_i$ any element of $Y_i$, but the only role of the elements of $H_i\setminus Y_i$
is to provide neighbours for the vertices in $U(\mathcal{T})$; so we can happily remove some of them if we also remove from 
$U(\mathcal{T})$  the vertices which no longer have neighbours in any of the (shrunken) sets $H_i$.

Let $B$ be the set of vertices in $U(\mathcal{T})$ with no neighbour in $Y(\mathcal{T})$. 
By \ref{privatecover}, since $Z(\mathcal{T})$ covers $B$, 
there exist $A'\subseteq Z(\mathcal{T})$ and $B'\subseteq B$ such that
\begin{itemize}
\item $A'$ covers $B\setminus B'$;
\item $B'$ is the union of $\delta\tau$ matching-covered sets;
\item every vertex in $B'$ has at most one neighbour in $A'$; and
\item every vertex in $A'$ has exactly $\delta\tau$ neighbours in $B'$.
\end{itemize}
Let $\Pi$ be the set of vertices in $B'$ that have a neighbour in $A'$; for $1\le i\le n$ let $H_i'=(H_i\cap A')\cup Y_i$;
and let $\mathcal{T}'$ be the template array with sequence $(Y_i, H_i')\;(1\le i\le n)$ and 
$U(\mathcal{T}')=(U(\mathcal{T})\setminus B')\cup \Pi$. Since $\chi(B')\le \delta\tau^2$ by (iii), this proves \ref{privatize}.~\bbox

The advantage of privatization is the following lemma, used when we have a shadowing of bounded degree.

\begin{thm}\label{commonroot}
Let $\eta\ge 1$ and $\zeta \ge \max(\eta,\alpha)$, and let $q,s\ge 0$.
Let $G$ satisfy {\rm (i)--(v)}, and 
let $\mathcal{T}$ be a 2-cleaned
$(\zeta,\eta)$-template array in $G$, with sequence $(Y_i,H_i)\;(1\le i\le n)$, that admits a privatization $\Pi$.
Let $B_1\ll B_n$ be a shadowing of degree at most $s$ relative to $U(\mathcal{T})\setminus \Pi$.
Let $1\le i\le n$, and let $\{D_j:j\in J\}$ be a bunch of daisies, each with root in $H_i$ and with $V(D_j)\cap \Pi=\emptyset$, with
$$|J|=2q\zeta\beta(1+(q+s)(\delta^2+1) + 2\delta+\delta\tau)\tau.$$
Then there exist $u\in H_i\cup B_i$ and $J'\subseteq J$ with $|J'|= q$, such that for each $j\in J'$,
$u$ is adjacent to the eye of $D_j$ and nonadjacent to the petals of $D_j$.
\end{thm}
\Proof
For each $j\in J$, let $u_j, v_j, P_j$ be the root, eye, and set of petals of $D_j$ respectively. 
Let $D$ be the digraph with vertex set $J_1$ in which for all distinct $j,j'\in J_3$, $j'$ is adjacent from $j$ if
$u_j$ is adjacent to $v_{j'}$. If some vertex of $D$ has outdegree at least $q$ we are done, so we assume not.
Hence by \ref{digraph} the graph underlying $D$ has chromatic number at most $2q$, and so there exists $J_1\subseteq J$
with $|J_1|=|J|/(2q)$ such that $u_j$ is nonadjacent to $v_{j'}$ for all distinct $j,j'\in J_1$.
In particular, the vertices $u_j(j\in J_1)$ are all distinct.

Since $|Y_i|=\zeta\eta$, it follows that $u_j\notin Y_i$ for at least $|J_1|-\zeta\eta$  values of $j\in J_1$.
Now each such $u_j$ has a neighbour in $Y_i$, and so there exist $J_2\subseteq J_1$
with 
$$|J_2|= (|J_1|-\zeta\eta)/|Y_i|= ((q+s)(\delta^2+1) + 2\delta+\delta\tau)\tau$$
and a vertex $y\in Y_i$, such that $u_j\notin Y_i$ and $y$ is adjacent to $u_j$ for each $j\in J_2$. 
Since $\{u_j:j\in J_1\}$ is $\tau$-colourable, there exists $J_3\subseteq J_2$ with 
$$|J_3|= |J_2|/\tau= (q+s)(\delta^2+1) + 2\delta+\delta\tau$$
such that the vertices $u_j(j\in J_3)$ are pairwise nonadjacent. Consequently for all distinct $j,j'\in J_3$
the daisies $D_j, D_{j'}$ are vertex-disjoint and no edge joins them. 

The set of vertices in $\Pi$ with distance two from $y$ has chromatic number at most $\tau$; fix some partition of this set
into $\tau$ stable sets.
For each $j\in J_3$, $u_j\in Z(\mathcal{T})$ and so $u_j$ has $\delta\tau$ neighbours in $\Pi$. 
We claim that all these neighbours belong to $B_i$.
For
$x\in \Pi$ be adjacent to $u_j$, and let $x\in B_k$ say. Then $x$ has a neighbour in $H_k$, from the definition of $B_k$;
but $x$ has a unique neighbour in $H(\mathcal{T})$, since $x\in \Pi$, and this neighbour in $u_j$, and so $u_j\in H_k$.
Since $u_j\in H_i$ it follows that $k=i$. Thus $u_j$ has $\delta\tau$ neighbours in $\Pi$, and they all belong to $B_i$.

The set of vertices in $\Pi$ with distance two from $y$ has chromatic number at most $\tau$; fix some partition of this set
into $\tau$ stable sets.
Consequently, for each $j\in J_3$, there are $\delta$ neighbours of $u_j$ that belong to the
same stable set of the partition. Since $|J_3|\ge \delta\tau$, there exists $J_4\subseteq J_3$ with 
$|J_4|= \delta$
and a stable subset $\Pi'$ of $\Pi\cap B_i$, such that for each $j\in J_4$, $u_j$ has at least $\delta$ neighbours in $\Pi'$.
For each $j\in J_4$, choose $\Pi_j\subseteq \Pi'\cap B_i$
with $|\Pi_j|=\delta$, such that every vertex in $\Pi_j$ is adjacent to $u_j$. Thus $G[\{y,u_j\}\cup \Pi_j]$ is a $(1,\delta)$-broom
with handle $y$, for each $j\in J_4$; and there are no edges between these brooms not incident with $y$.

Let $\{y\}\cup \bigcup_{j\in J_4}\Pi_j=Q$ say.
If some vertex in $Q$ is adjacent to $q$ of the vertices $v_j(j\in J_3)$ we are done, so we assume not.
Also by hypothesis each vertex in $Q$ has a neighbour in $P_{j}$ for at most $s$ values of $j\in \{1\ll n\}$.
Since $|Q|=\delta^2+1$, there are at most $(q+s)(\delta^2+1)$ values of $j\in J_3\setminus J_4$ such that some vertex in $Q$
has a neighbour in $P_j\cup \{v_j\}$. Since $|J_3|-|J_4|-(q+s)(\delta^2+1)\ge \delta$, 
there exists $J_5\subseteq J_3\setminus J_4$
with 
$|J_5|= \delta$
such that for each $j\in J_5$, no vertex in $Q$ is adjacent to $v_j$ or has a neighbour in $B_j\setminus \Pi$.
It follows that 
$G[\{y,u_j,v_j\}\cup P_j]$ is a $(2,\delta)$-broom with handle $y$ for each $j\in j_5$. By taking the union of the $(1,\delta)$-brooms
$G[\{y,u_j\}\cup \Pi_j]$ for each $j\in J_4$ and the $(2,\delta)$-brooms $G{\{y,u_j,v_j\}\cup P_j}]$ for each $j\in J_5$, 
we find that $G$ contains $T(\delta)$, a contradiction. This proves \ref{commonroot}~\bbox

\section{Edges between $H(\mathcal{T})$ and $U(\mathcal{T})$.}

Our next goal is to bound the number of neighbours each vertex of $U(\mathcal{T})$ has in $H(\mathcal{T})$.

\begin{thm}\label{manyHvert}
Let $\eta\ge \delta$ and $\zeta\ge \max(\eta,\alpha)+\delta$ be integers; then there exists $\ell\ge 0$ with the following
property. Let $G$ satisfy {\rm (i)--(v)}, and let $\mathcal{T}$ be a 2-cleaned
$(\zeta,\eta)$-template array in $G$, with sequence $(Y_i,H_i)\;(1\le i\le n)$, that admits a privatization.
Let $X$ be the set of vertices in $U(\mathcal{T})$ that have at least $(\beta+1)\gamma\delta$
neighbours in $H(\mathcal{T})$.
Then
$\chi(X)\le \ell$.
\end{thm}
\Proof
Let
\begin{eqnarray*}
s&=& 2(2\gamma+1)\delta\tau+\gamma,\\
q&=&2\delta (2\gamma(\delta+1)+1)+\gamma,\\ 
m&=&2q\zeta\beta(1+(q+s)(\delta^2+1) + 2\delta+\delta\tau)\tau, \text{ and}\\
\ell&=&2sm\zeta^2\delta(2(\delta+2)s+3)\beta^2\tau^3. 
\end{eqnarray*}
We claim that $\ell$ satisfies the theorem. For let $\mathcal{T}$, 
$(Y_i,H_i)\;(1\le i\le n)$ and $X$ be as in the theorem.
For each $v\in X$, since $v$ has at least $(\beta+1)\gamma\delta$ neighbours in $H(\mathcal{T})$,
and there are only at most $\gamma$ values of $i$ such that $v$ has a neighbour in $H_i$, it follows that there exists $i\in \{1\ll n\}$
such that $v$ has at least $(\beta+1)\delta$ neighbours in $H_i$. 
Choose a shadowing $B_1\ll B_n$ such that for $1\le i\le n$, every vertex in 
$B_i\cap X$ has at least $(\beta+1)\delta$ neighbours in $H_i$.
\\
\\
(1) {\em The shadowing $B_1\ll B_n$ has degree less than $s$.}
\\
\\
Suppose not, and choose $y\in V(\mathcal{T})$
and $J\subseteq \{1\ll n\}$ with $|J|=s$ such that for each $j\in J$ there exists $u_j\in X\cap B_j$
adjacent to $y$. Since $y$ has a neighbour in $H_j$ for at most $\gamma$ values of $j$,
there exists $J_1\subseteq J$ with $|J_1|=|J|-\gamma$ such that $y$ has no neighbour in $H_j$ for each $j\in J_1$.
Since the subgraph induced on $\{u_j:j\in J_1\}$ is $\tau$-colourable, there exists $J_2\subseteq J_1$
with $|J_2|=|J_1|/\tau=2(2\gamma+1)\delta$ such that the vertices $u_j(j\in j_1)$ are pairwise nonadjacent.
Let $D$ be the digraph with vertex set $J_2$ in which for all distinct $j,j'\in J_2$, $j'$ is adjacent from $j$
if $u_j$ has a neighbour in $H_{j'}$. Since $D$ has maximum outdegree at most $\gamma$, by \ref{digraph} the graph underlying $D$
is $(2\gamma+1)$-colourable, and so there exists $J_3\subseteq J_2$ with $|J_3|=|J_2|/(2\gamma+1)=2\delta$
such that for all distinct $j,j'\in J_2$, $u_j$ has no neighbour in $H_{j'}$. 
For each $j\in J_3$, since $u_j\in X$, $u_j$ has at least $\delta(\beta+1)$ neighbours in $H_j$; 
and since $H_j$ is $(\beta+1)$-colourable, there is a stable set
$P_j$ of $\delta$ such neighbours. Thus $G[\{y,u_j\}\cup P_j]$ is a $(1,\delta)$-broom with handle $y$.
Let $v_j$ be a neighbour of $u_j$ in $H_j$, choosing $v_j\in Y_j$ if possible. If $v_j\in Y_j$,
let $A$ be a part of $Y_j$ not containing $v_j$; then since $u_j\notin H_j$, $u_j$ has at most $\eta-1$ neighbours in $A$, and 
since $\zeta\ge \eta-1 +\delta$, there is a set $Q_j\subseteq A$
of $\delta$ vertices all nonadjacent to $u_j$. If $v_j\notin Y_j$, then $u_j$ has no neighbour in $Y_j$, and 
since $v_j$ is $\eta$-mixed on $Y_j$ and $h\ge \delta$, it follows that $v_j$
has a set $Q_j$ of $\delta$ neighbours in some part of $Y_j$. In either case $G[\{y,u_j,v_j\}\cup Q_j]$ is a $(2,\delta)$-broom
with handle $y$. By choosing the $(1,\delta)$-broom for $\delta$ values of $j\in J_3$, and the $(2,\delta)$-broom
for the other $\delta$ values of $j\in J_3$, and taking their union, we find that $G$ contains $T(\delta)$, a contradiction.
This proves (1).
\\
\\
(2) {\em There is no bunch of daisies $\{D_j:j\in J\}$ with $|J|=m$ such that $V(D_j)\subseteq X\cup H(\mathcal{T})$
for each $j\in J$.}
\\
\\
Suppose such a bunch $\{D_j:j\in J\}$ exists. Let $i\in \{1\ll n\}$ such that the root of $D_j$ belongs to $H_i$
for each $j\in J$. Let $\Pi$ be a privatization. (Thus $\Pi\cap X=\emptyset$, since $(\beta+1)\gamma\delta\ge 2$.)
By \ref{commonroot} applied to the $(\zeta,\eta)$-template array $\mathcal{T}$ with sequence 
$(Y_j,H_j)(1\le j\le n)$
and $U(\mathcal{T}') = X\cup \Pi$, there exist $y\in H_i\cup B_i$ and $J_1\subseteq J$ with $|J_1|= q$, such that for each $j\in J_1$,
$y$ is adjacent to the eye ($u_j$ say) of $D_j$ and has no neighbour in the set of petals ($P_j$ say) of $D_j$. 
Thus $G[\{y,u_j\}\cup P_j]$ is a $(1,\delta)$-broom with handle $y$, for each $j\in J_1$.
Since there are at most $\gamma$ values of $j\in J_1$ such that $y$ has a neighbour in $H_j$, there exists $J_2\subseteq J_1$
with $|J_2|=|J_1|-\gamma$ such that $y$ has no neighbour in $H_j$ for $j\in J_2$. 

Let $D$ be the digraph with vertex set $J_2$ in which for all distinct $j,j'\in J_2$, $j'$ is adjacent from $j$
if some vertex in $\{u_j\}\cup P_j$ has a neighbour in $H_{j'}$. Since $D$ has maximum outdegree at most $\gamma(\delta+1)$, 
by \ref{digraph} the graph underlying $D$
is $(2\gamma(\delta+1)+1)$-colourable, and so there exists $J_3\subseteq J_2$ with $|J_3|=|J_2|/(2\gamma(\delta+1)+1)=2\delta$
such that for all distinct $j,j'\in J_3$, no vertex in $\{u_j\}\cup P_j$ has a neighbour in $H_{j'}$. 

For each $j\in J_3$, choose a neighbour $v_j$ of $u_j$, such that 
\begin{itemize}
\item if $u_j$ has a neighbour in $Y_j$ then $v_j\in Y_j$;
\item if $u_j$ has no neighbour in $Y_j$ and has a neighbour in $H_j$ then $v_j\in H_j$;
\item if $u_j$ has no neighbour in $H_j$ then $v_j\in P_j$.
\end{itemize}
We claim that in each case, $y$ is nonadjacent to $v_j$, and there is a stable set $Q_j\subseteq H_j$ of neighbours of $v_j$,
all nonadjacent to $u_j$, with $|Q_j|=\delta$. To see this, if $v_j\in H_j$ the proof is as in the proof of (1),
so we assume that $v_j\in P_j$, and therefore $u_j$ has no neighbours in $H_j$. Since $P_j\subseteq X$,
$v_j$ has at least $\delta(\beta+1)$ neighbours in $H_j$, and since $H_j$ is $(\beta+1)$-colourable, the claim follows.
In particular, $G[\{y,u_j,v_j\}\cup Q_j]$ is a $(2,\delta)$-broom for each $j\in J_3$.
By choosing the $(1,\delta)$-broom for $\delta$ values of $j\in J_3$, and the $(2,\delta)$-broom
for the other $\delta$ values of $j\in J_3$, and taking their union, we find that $G$ contains $T(\delta)$, a contradiction.
This proves (2).

\bigskip
From (1), (2) and \ref{getdaisies}, it follows that
$\chi(X)\le  \ell$.
This proves \ref{manyHvert}.~\bbox

\bigskip
The bound $(\beta+1)\gamma\delta$ will be very useful in the remainder of the proof, and for convenience let us
define $\epsilon=(\beta+1)\gamma\delta$, for the remainder of the paper.
Let us say a 2-cleaned $(\zeta,\eta)$-template array $\mathcal{T}$ in $G$ is {\em 3-cleaned} if
every vertex in $U(\mathcal{T})$ has fewer than $\epsilon$ neighbours in $H(\mathcal{T})$.
We deduce:

\begin{thm}\label{get3cleaned}
Let $\eta\ge \delta$ and $\zeta\ge \max(\eta,\alpha+\delta)$ be integers.
Then there is a non-decreasing
function $\phi:\mathbb{N}\rightarrow \mathbb{N}$ with the following property.
For all $c\ge 0$, if $G$ satisfies {\rm (i)--(v)} and $\chi(G)>\phi(c)$ then there is a
$3$-cleaned $(\zeta,\eta)$-template array $\mathcal{T}$
in $G$ such that $\chi(U(\mathcal{T}))>c$.
\end{thm}
\Proof 
Let $\psi$ satisfy \ref{get2cleaned} (with $\phi$ replaced by $\psi$), and let $\ell$ be as in \ref{manyHvert}.
For all $c\ge 0$ define $\phi(c)=\psi(c+\delta\tau^2+\ell)$;
we claim this satisfies the theorem. For let $c\ge 0$, and let $G$ satisfy {\rm (i)--(v)}
with $\chi(G)>\phi(c)$. By \ref{get2cleaned},
there is a
$2$-cleaned $(\zeta,\eta)$-template array $\mathcal{T}$
in $G$ such that $\chi(U(\mathcal{T}))>c+\delta\tau^2+\ell$.
By \ref{privatize}, there
is a 2-cleaned $(\zeta,\eta)$-template array $\mathcal{T}_1$ in $G$
that admits a privatization, such that $\chi(U(\mathcal{T}_1))> c+\ell$.
Let $X$ be the set of vertices in $U(\mathcal{T}_1)$ that have at least $\epsilon$
neighbours in $H(\mathcal{T}_1)$. 
Let $\ell$ be as in \ref{manyHvert}; then by \ref{manyHvert}, 
$\chi(X)\le \ell$. Let $\mathcal{T}'$ be the $(\zeta,\eta)$-template array with the same sequence as $\mathcal{T}_1$
and with $U(\mathcal{T}') = U(\mathcal{T}_1)\setminus X$. 
It follows that $\mathcal{T}'$ is $3$-cleaned, and
$\chi(U(\mathcal{T}'))\ge \chi(U(\mathcal{T}_1))-\ell> c.$ 
This proves \ref{get3cleaned}.~\bbox

\section{Edges within a shadowing}

Now we have come to the final stage of the proof: we investigate the edges between different sets of a shadowing.
First we prove that there is some template array such that every shadowing has bounded degree; and privatize it;
and then we will show that
for a privatized template array, if every shadowing has bounded degree then the graph has bounded chromatic number.

So far, our technique in this paper has been to start with a template array, and make nicer and nicer ones at the cost
of reducing the chromatic number of $U(\mathcal{T})$. This has more-or-less reached its limit, with \ref{get3cleaned},
so now we need to do something different. To prove the next result, we will start with a 3-cleaned
template array $\mathcal{T}$, and apply \ref{get3cleaned} to $G[U(\mathcal{T})]$  to get a second one, with vertex set
a subset of $U(\mathcal{T})$; and repeat, generating a nested sequence of template arrays.

\begin{thm}\label{nested}
Let $\eta\ge \alpha+ 2(\delta+1)^3(\epsilon+1)^2$ and $\zeta\ge \eta+\delta$ be integers.
Then there is a non-decreasing
function $\phi:\mathbb{N}\rightarrow \mathbb{N}$, and an integer $s$, with the following property.
Let $G$ satisfy {\rm (i)--(v)}, with  $\chi(G)> \phi(c)$.
Then there is a $3$-cleaned $(\zeta,\eta)$-template array $\mathcal{T}$
in $G$, with sequence $(Y_i,H_i)\;(1\le i\le n)$, such that $\chi(U(\mathcal{T}))>c$, and 
such that for each vertex $v\in U(\mathcal{T})$, there are fewer than
$s$ values of $i\in \{1\ll n\}$ such that some neighbour of $v$ in $U(\mathcal{T})$ has a neighbour in $H_i$.
\end{thm}
\Proof

Let 
\begin{eqnarray*}
s_3 &=& (\delta(\delta+1)+1)\epsilon +\delta\\ 
s_2 &=& ((2\delta\epsilon+1)\delta\tau+\epsilon)s_3\\
s_1&=& (2\epsilon+1)\tau s_2\\
s&=&s_1+\epsilon\\
t_4&=&2\delta\\
t_3& =& 2\delta t_4\\
t_2& =& (\alpha\tau2^{\beta\zeta}+1)t_3\\
t_1&=&t_2^{s_2}, \text{ and} \\
t&=&1+2^{s_2}\delta\tau+t_1.
\end{eqnarray*}
Let $\psi$ satisfy \ref{get3cleaned} (with $\phi$ replaced by $\psi$). Define $\phi^0(x)=x$ for $x\ge 0$, and inductively for
$i\ge 1$, let $\phi^i(x) = \psi(\phi^{i-1}(x))$ for $x\ge 0$. Let $\phi = \phi^t$; 
 we claim that $\phi,s$ satisfy
the theorem. 

Let $G$ satisfy {\rm (i)--(v)}, with  $\chi(G)> \phi(c)$. Define $U^0=V(G)$;
thus $\chi(U_0)>\phi^{t}(c)$. For $1\le j\le t$
define $\mathcal{T}^j$ and $U^j$ inductively as follows. Let $j$ satisfy $1\le j\le t$, and suppose we have defined $U^{j-1}$, and
$\chi(U^{j-1})>\phi^{t-j+1}(c)$. By \ref{get3cleaned} applied to $G[U^{j-1}]$, 
there is a
$3$-cleaned $(\zeta,\eta)$-template array $\mathcal{T}^j$
in $G[U^{j-1}]$ such that $\chi(U(\mathcal{T}^j))>\phi^{t-i}(c)$. Let $U^j= U(\mathcal{T}^j)$. This completes the inductive definition.
Let $\mathcal{T}$ be the $(\zeta,\eta)$-template array with the same sequence as $\mathcal{T}^1$, and with $U(\mathcal{T}) = U^t$;
we will show that $\mathcal{T}$ satisfies the theorem. Certainly it is 3-cleaned, and $\chi(U(\mathcal{T}))>c$.

We remark that for all $j<j'\le t$, every vertex in $H(\mathcal{T}^{j'})$ has fewer than $\epsilon$ neighbours in $H(\mathcal{T}^j)$
(since $H(\mathcal{T}^{j'})\subseteq U(\mathcal{T}^j)$ and $\mathcal{T}^j$ is 3-cleaned),
and for all $j$ 
every vertex in $U^t$
has fewer than $\epsilon$ neighbours in $H(\mathcal{T}^j)$ (for the same reason); but when $j'>j$ we know nothing about the number of 
neighbours a vertex in $H(\mathcal{T}^j)$
has in $H(\mathcal{T}^{j'})$.

For $1\le j\le t$, let the sequence of $\mathcal{T}^j$ be $(Y_i^j,H_i^j)\;(1\le i\le n_j)$. We assume for a contradiction 
that there exist $y\in U^t$ and a subset $I\subseteq \{1\ll n_1\}$ with $|I|=s$, such that for each $i\in I$,
there exists $u_i\in U^t$ adjacent to $y$ and with a neighbour in $H_i^1$. Since there are at most $\epsilon$
values of $i$ such that $y$ has a neighbour in $H_i^1$, there exists $I_1\subseteq I$
with $|I_1|=|I|-\epsilon=s_1$ such that for each $i\in I_1$, $y$ has no neighbour in $H_i^1$.
\\
\\
(1) {\em There exist $I_2\subseteq I_1$ with $|I_2|=s_2$, and for each $i\in I_2$, a $(2,\delta)$-broom $B_i$ of 
$G[\{y,u_i\}\cup H_i^1]$ with handle $y$, such that every edge joining two of these brooms is incident with $y$.}
\\
\\
Let $D$ be the digraph with vertex set $I_1$
in which all for all distinct $i,i'\in I_1$, $i'$ is adjacent from $i$ if $u_i$ has a neighbour in $H_{i'}^1$. 
Since $D$ has maximum outdegree
at most $\epsilon$, by \ref{digraph} the graph underlying $D$ has chromatic number at most $2\epsilon+1$. Hence there exists
$I_1'\subseteq I_1$ with $|I_1'|=s_1/(2\epsilon+1)$, such that for all all distinct $i,i'\in I_1'$,
$u_i$ has no neighbour in $H_{i'}^1$. Since the set $\{u_i:i\in I_1'\}$ has chromatic number at most $\tau$ (because each $u_i$
is adjacent to $y$), there exists
$I_2\subseteq I_1'$ with $|I_2|=|I_1'|/\tau$ such that the vertices $u_i(i\in I_2)$ are pairwise nonadjacent.
For each $i\in I_2$, choose $w_i\in H_i^1$ adjacent to $u_i$, with $w_i\in Y_i^1$ if possible. If $w_i\in Y_i^1$,
let $A$ be a part of $Y_i^1$ not containing $w_i$; then since $u_i$ has fewer than $\eta$ neighbours in $A$, and $w_i$ is adjacent
to every vertex in $A$, there exists $R_i\subseteq A$ of cardinality $\delta$ all nonadjacent to $u_i$. If $w_i\notin Y_i^1$,
then $u_i$ has no neighbour in $Y_i$; since $w_i\in H_i^1$, $w_i$ has at least $\eta\ge \delta$ neighbours in some part of $Y_i^1$;
and so there exists a stable set $R_i\subseteq Y_i^1$ of cardinality $\delta$, 
all adjacent to $w_i$ and not to $u_i$. In either case,
$G[\{y,u_i,w_i\}\cup R_i]$ is a $(2,\delta)$-broom with handle $y$. This proves (1).

\bigskip

Let $J_1$ be the set of all $j\in \{2\ll t\}$ such that some vertex in $H(\mathcal{T}^j)$ is adjacent to at least $s_3$
of the vertices $u_i(i\in I_2)$. 
\\
\\
(2) {\em If $j\in \{2\ll t\}\setminus J_1$, there is a subset $X^j\subseteq H(\mathcal{T}^j)$ and a subset
$I^j\subseteq I_2$, such that 
\begin{itemize}
\item $|X^j|=|I^j|=(2\delta\epsilon+1)\delta\tau$; 
\item $y$ has no neighbour in $X^j$;
\item every vertex in $X^j$ has a unique neighbour in $\{u_i:i\in I^j\}$; and
\item every vertex in $\{u_i:i\in I^j\}$ has a unique neighbour in $X^j$.
\end{itemize}
}
Since for all $i\in I_2$, $u_i$ has a neighbour in $H(\mathcal{T}^j)$, there exists $X\subseteq H(\mathcal{T}^j)$,
minimal such that each $u_i(i\in I_2)$ has a neighbour in $X$. Since every vertex in $X$ is adjacent to fewer than $s_3$
vertices in $\{u_i:i\in I_2\}$, it follows that $|X|\ge s_2/s_3=(2\delta\epsilon+1)\delta\tau+\epsilon$. 
From the minimality of $X$, for each $x\in X$
there exists $i(x)\in I_2$ such that $x$ is the unique neighbour of $u_{i(x)}$ in $X$. Since $y$
has at most $\epsilon$ neighbours in $H(\mathcal{T}^j)$, 
there exists $X^j\subseteq X$ with 
$|X^j|=(2\delta\epsilon+1)\delta\tau$ such that $y$ has no neighbours in $X^j$. 
Let $I^j=\{i(x):x\in X^j\}$; then this proves (2).
\\
\\
(3) {\em $|\{2\ll t\}\setminus J_1|\le 2^{s_2}\delta\tau$.}
\\
\\
Suppose not. For each $j\in \{2\ll t\}\setminus J_1$, there are at most $2^{s_2}$ possibilities for the set $I^j$, 
and so there exists 
$J'\subseteq  \{2\ll t\}\setminus J_1$ with $|J'|= \delta\tau$, and a subset $I_3\subseteq I_2$ (necessarily 
with $|I_3|=(2\delta\epsilon+1)\delta\tau$),
such that $I^j=I_3$ for each $j\in J'$. Let $X=\bigcup_{j\in J'}X^j$; then $\chi(X)\le \tau$.
Take a partition $Z_1\ll Z_{\tau}$ of $X$ into stable sets. For each $i\in I_3$, since $u_i$ has $\delta\tau$ neighbours in
$X$, there exists $r\in \{1\ll \tau\}$ such that $u_i$ has at least $\delta$ neighbours in $Z_r$. Since there
are only $\tau$ possibilities for $r$, there exist $I_4\subseteq I_3$ with $|I_4|= |I_3|/\tau=(2\delta\epsilon+1)\delta$,
and $r\in \{1\ll \tau\}$, such that for each $i\in I_4$, $u_i$ has at least $\delta$ neighbours in $Z_r$; let $N_i$
be a set of $\delta$ such neighbours.
Since for each $j\in J'$, every vertex in $X_j$ has a unique neighbour in $\{u_i:i\in I_2\}$, the sets $N_i(i\in I_4)$
are pairwise disjoint and their union is stable. Let $D$ be the digraph with vertex set $I_4$, in which for all distinct $i,i'\in I_4$,
$i'$ is adjacent from $i$ if some vertex in $N_i$ has a neighbour in $V(B_{i'})$ (as defined in (1)). 
Then $D$ has maximum outdegree
$\delta\epsilon$, since $|N_i|=\delta$ and each vertex in $N_i$ has at most $\epsilon$ neighbours in $H^1$
(and none in $V(B_{i'})\setminus H^1$). By \ref{digraph} the graph underlying $D$ has chromatic number at most
$2\delta\epsilon+1$, and so there exists $I_5\subseteq I_4$ with 
$$|I_5|= |I_4|/(2\delta\epsilon+1)=\delta,$$ 
such that for all distinct $i,i'\in I_5$, no vertex in $N_i$ has a neighbour in $V(B_{i'})$.
It follows that $G$ contains $T(\delta)$,
a contradiction. This proves (3).

\bigskip
From (3), it follows that $|J_1|\ge t-1-2^{s_2}\delta\tau=t_1$. For each $j\in J_1$, choose $v^j\in H(\mathcal{T}^j)$
adjacent to at least $s_3$  of the vertices $u_i(i\in I_2)$. Consequently there exists $J_2\subseteq J_1$
with $|J_2|=t_2\le |J_1|2^{-s_2}$ and a subset $I_3\subseteq I_2$ with $|I_3|=s_3$, such that every vertex $v^j(j\in J_2)$
is adjacent to every vertex $u_i(i\in I_3)$. 
For each $j\in J_2$, there exists $i\in \{1\ll n_j\}$ such that $v^j\in H_i^j$; let us write $H^j=H_i^j$ and $Y^j=Y_i^j$,
since we will have no need for the other terms of the sequence of $\mathcal{T}^j$.
\\
\\
(4) {\em There exists $J_3\subseteq J_2$ with $|J_3| =|J_2|/(\alpha\tau2^{\beta\zeta}+1)=t_3$ 
such that for all distinct $j,j'\in J_2$,
$v_j$ is not dense to $Y^{j'}$.}
\\
\\
Let $D$ be the digraph with vertex set $J_1$ in which for all distinct $j,j'\in J_1$, $j'$ is adjacent from $j$ if 
$v_j$ is dense to $Y^{j'}$. Since $v_j$ is not dense to $Y^{j'}$ if $j>j'$, it follows that
$D$ is acyclic, and has maximum indegree at most $\alpha\tau2^{\beta\zeta}$, by \ref{dense},
and so by \ref{digraph} the graph underlying $D$ has chromatic number at most $\alpha\tau2^{\beta\zeta}+1$. This proves (4).
\\
\\
(5) {\em For each $j_0\in J_3$ there are fewer than $\delta$ values of $j\in J_3\setminus \{j_0\}$ such that
$v_{j_0}$ has at least $\delta(\delta+1)\epsilon+s_3\epsilon$ neighbours in some part of
$Y^{j}$.}
\\
\\
For suppose that there exists $J_4\subseteq J_3\setminus \{j\}$ with $|J_4|=\delta$ such that for each $j\in J_4$,
$v_{j_0}$ has at least $\delta(\delta+1)\epsilon+s_3\epsilon$ neighbours in some part of $Y^{j}$. 
For each $j\in J_4$, since there is a part of $Y^{j}$ in which
$v_{j_0}$ has fewer than $\alpha$ neighbours (because $v_{j_0}$ is not dense to $Y^{j}$), 
it follows that there are distinct parts $A,A'$
of $Y^{j}$, such that $v_{j_0}$ has at least $\delta(\delta+1)\epsilon+s_3\epsilon$ neighbours in $A$ and at least $\zeta-\alpha+1\ge \delta(\delta+1)\epsilon+s_3\epsilon$
non-neighbours in $A'$. Since at most $s_3\epsilon$
vertices of $Y^j$ have a neighbour in $\{u_i:i\in I_3\}$, there is a subset $P^j\subseteq A$ with cardinality
$\delta(\delta+1)\epsilon$, such that all vertices in $P^j$ are adjacent to $v_{j_0}$ and have no neighbours in
$\{u_i:i\in I_3\}$; and there is a subset $Q^j\subseteq A'$ with cardinality 
$\delta(\delta+1)\epsilon$, such that all vertices in $Q^j$ are nonadjacent to $v_{j_0}$ and have no neighbours in 
$\{u_i:i\in I_3\}$.

For each $j\in J_4$, we choose a $(1,\delta)$-broom $C^{j}$ of $G[Y^{j}\cup \{v_{j_0}\}]$ with handle $v_{j_0}$, 
inductively as follows. Let $j\in J_4$,
and assume that $C^{j'}$ is defined for all $j'\in J_4$ with $j'>j$.
Let $S$ be the union of all the sets $V(C^{j'})\setminus \{v_{j_0}\}$ for $j'\in J_4$ with $j'>j$.
Then $|S|\le  (\delta-1)(\delta+1)$,
and since each vertex in $S$ has at most $\epsilon$ neighbours in $H^j$, it follows that at most
$(\delta-1)(\delta+1)\epsilon$ vertices in $H^j$ have a neighbour in $S$. Since
$$|P^j|,|Q^j|=\delta(\delta+1)\epsilon> (\delta-1)(\delta+1)\epsilon+\delta$$
there exist a vertex in $P^j$, and a set of $\delta$ vertices in $Q^j$, each with no neighbours in $S$. Hence
there is a $(1,\delta)$-broom $C^{j}$ of $G[Y^{j}\cup \{v_{j_0}\}]$ containing no neighbours of $S$ different from $v_{j_0}$.
This completes the inductive definition of $C^j$ for $j\in J_4$.

Now let $S$ be the union of the sets $V(C^{j})$ for $j\in J_4$. Then $|S|=\delta(\delta+1)+1$,
and so there are at most $(\delta(\delta+1)+1)\epsilon$ vertices in $H^1$ with neighbours in $S$. 
Since $s_3\ge (\delta(\delta+1)+1)\epsilon +\delta$, there are at least $\delta$ values of $i\in I_3$
such that the edge $u_iv_{j_0}$ is the only edge between $V(B_i)$ and $V(C^j)$ for each $j\in J_4$; and so $G$
contains $T(\delta)$ (with handle $v_{j_0}$), a contradiction. This proves (5).
\\
\\
(6) {\em There exists $J_4\subseteq J_3$ with $|J_4|=t_4$ such that for all distinct $j,j'\in J_4$,
$v_{j}$ has fewer than $\delta(\delta+1)\epsilon+s_3\epsilon$ neighbours in each part of $Y^{j}$.}
\\
\\
Let $D$ be the digraph with vertex set $J_3$ in which for all distinct $j,j'\in J_3$, $j'$ is adjacent from $j$ if 
$v_j$ has at least $\delta(\delta+1)\epsilon+s_3\epsilon$ neighbours in some part of $Y^{j}$. By (5),
$D$ has maximum outdegree at most $\delta-1$, and so by \ref{digraph} the graph underlying $D$ has chromatic number at most  $2\delta$,
and the claim follows. This proves (6).

\bigskip
Fix $i\in I_3$; and partition $J_4$ into two sets $J_4', J_4''$, both of cardinality $\delta$.
For each $j\in J_4'$, we define a $(2,\delta)$-broom $C^j$ of $G[Y^j\cup \{v_j,u_i\}]$ with handle $u_i$,
and for each $j\in J_4''$, we define a $(1,\delta)$-broom $C^j$ of $G[Y^j\cup \{v_j,u_i\}]$ with handle $u_i$,
inductively as follows.
Let $j\in J_4$,
and assume that $C^{j'}$ is defined for all $j'\in J_4$ with $j'>j$.
Let $S$ be the union of $\{u_i\}$ and all the sets $V(C^{j'})\setminus \{v^{j'}\}$ for $j'\in J_4$ with $j'>j$.
Then $|S|\le  (\delta+1)(t_4-1)$.
There are distinct parts $A,A'$ of $Y^j$ such that $v_j$ has at least $\eta$ neighbours in $A$
and at most $\alpha-1$ neighbours in $A'$. But
at most $\epsilon(\delta+1)(t_4-1)$ vertices of $A$ have a neighbour in $S$, and at most 
$(t_4-1)(\delta(\delta+1)\epsilon+s_3\epsilon)$ vertices in $A$ have a neighbour
in $\{v_{j'}:j'\in J_4\setminus \{j\}\}$, by (6) and the definition of $J_4$, and the same for $A'$.
Since
$$\eta, \zeta-\alpha+1\ge  (\delta+1)(t_4-1)\epsilon +(t_4-1)(\delta(\delta+1)\epsilon+s_3\epsilon)+ \delta,$$
there exist a set of $\delta$ neighbours of $v_j$ in $A$, and a set of $\delta$ non-neighbours of $v_j$ in $A'$, 
all with no neighbours in $S\cup \{v_{j'}:j'\in J_4\setminus \{j\}\}$.
Consequently the desired broom $C_j$ can be chosen as specified. This completes the inductive definition.
But by taking the union of all the $C^j(j\in J_4)$ , we see that $G$ contains $T(\delta)$, a contradiction. This
proves \ref{nested}.~\bbox

\begin{thm}\label{manyB}
Let $\zeta,\eta,s\ge 0$.
Let $\mathcal{T}$ be a 2-cleaned
$(\zeta,\eta)$-template array in $G$, with sequence $(Y_i,H_i)\;(1\le i\le n)$, that admits a privatization $\Pi$.
Suppose that for each vertex $v\in U(\mathcal{T})$, there are fewer than
$s$ values of $i\in \{1\ll n\}$ such that some neighbour of $v$ in $U(\mathcal{T})$ has a neighbour in $H_i$.
Let $q=2\delta+s$ and
$$r=(4(\delta+1)s+1)q\zeta\beta\tau(1+\tau((q+s)\delta^2+(2s(\delta+1)+1)\delta\tau)).$$
Then
$$\chi(U(\mathcal{T})\setminus \Pi)\le 3rs\beta\delta\zeta\tau^2.$$
\end{thm}

\Proof It follows that every shadowing has degree less than $s$.
For $1\le i\le n$, let $B_i$ be the set of vertices in $U(\mathcal{T})$ with a neighbour in $H_i$ and with no neighbour
in $H_1\cup\cdots\cup H_{i-1}$. It follows that $(B_1\ll B_n)$ is a shadowing, and hence has degree less than $s$.

\bigskip
For distinct $u,v\in U(\mathcal{T})$, we say that $v$ is {\em later} than $u$ if $u\in B_i$ and $v\in B_j$ where $j>i$.
For $i,a,b,c\in \{1\ll n\}$, we say that $i$ is {\em strong} to $(a,b,c)$ if 
\begin{itemize}
\item $i<\min(a,b,c)$ and $a<c$ ;
\item there exist $u\in B_i\setminus \Pi$ and $v\in B_a\setminus \Pi$, adjacent;
\item there exist $\delta$ vertices in $B_b\setminus \Pi$, pairwise nonadjacent, and all adjacent to $u$ 
(and possibly also adjacent to $v$); and 
\item there exist $\delta$ vertices in $B_c\setminus \Pi$, pairwise nonadjacent, and all adjacent to $v$ and not to $u$. 
\end{itemize}
(1) {\em For $1\le i\le n$, there do not exist $r$ triples $(a_1,b_1,c_1)\ll (a_r,b_r,c_r)$, such that $i$ is strong to them all, and
$a_j,b_j,c_j$ are different from $a_{j'}, b_{j'}, c_{j'}$ for $1\le j<j'\le r$.}
\\
\\
Suppose that $r$  such triples exist. For $1\le j\le r$, choose $u_j\in B_i\setminus \Pi$ and $v_j\in B_{a_j}\setminus \Pi$ adjacent to $u_j$,
and a stable set $P_j$ of $\delta$ vertices in $B_{b_j}\setminus \Pi$, all adjacent to $u_j$, and a stable set $Q_j$ 
of $\delta$ vertices in $B_{r_j}\setminus \Pi$, all adjacent to $v_j$ and not to $u_j$. Let $D$ be the digraph
with vertex set $\{1\ll r\}$ in which for all distinct $j,j'\in \{1\ll r\}$, $j'$ is adjacent from $j$ if 
some vertex $u \in \{u_j,v_j\}\cup P_j\cup Q_j$ is adjacent in $G$ to a vertex $v\in \{v_{j'}\}\cup P_{j'}\cup Q_{j'}$
and $u$ is earlier than $v$.
Then $D$ has maximum outdegree less than 
$2(\delta+1)s$, and so by \ref{digraph} there exists $J\subseteq\{1\ll r\}$ with 
$$|J|\ge r/(4(\delta+1)s)= q\zeta\beta(1+\tau((q+s)\delta^2+(2s(\delta+1)+1)\delta\tau))$$
such that for all distinct $j,j'\in J$, there are no edges between $\{u_j,v_j\}\cup P_j\cup Q_j$ and 
$\{v_{j'}\}\cup P_{j'}\cup Q_{j'}$. In particular, the vertices $u_j(j\in J)$ are all distinct.
Since the set $\{u_j:j\in J\}$ has chromatic number at most $\tau$, there exists $J_1\subseteq J$
with $|J_1|= |J|/\tau$ such that the vertices $u_j(j\in J_1)$ are pairwise nonadjacent.
For each $j\in J_1$ choose $w_j\in H_j$ adjacent to $u_j$; then $w_j$ has no neighbours in $P_j$, from the definition of
the shadowing, and so $G[\{w_j,u_j\}\cup P_j]$ is a daisy $D_{b_j}$ say, and $\{D_{b_j}:j\in J_1\}$ is a bunch of daisies.
But 
$$|J|/\tau\ge 2q\zeta\beta(1+(q+s)(\delta^2+1) + 2\delta+\delta\tau)\tau,$$
and so by \ref{commonroot}, applied to $\{D_{b_j}:j\in J_1\}$, we deduce that
there exist $w\in H_i\cup B_i$ and $J_2\subseteq J_1$ with $|J_2|= q$, such that for each $j\in J_2$,
$w$ is adjacent to the eye $u_j$ of $D_{b_j}$ and nonadjacent to the petals $P_j$ of $D_{b_j}$.
Moreover there are no edges between $\{u_j,v_j\}\cup P_j\cup Q_j$ and
$\{u_{j'},v_{j'}\}\cup P_{j'}\cup Q_{j'}$ for all distinct $j,j'\in J_2$. Now $w$ has neighbours in at most $s$
of the sets $\{v_j\}\cup P_j\cup Q_j\;(j\in J_2)$, and so there exists $J_3\subseteq J_2$ with $|J_3|=2\delta$ such that $w$
has no neighbours in $\{v_j\}\cup P_j\cup Q_j$ for $j\in J_3$. Hence for $j\in J_3$, $G[\{w,u_j\}\cup P_j]$ is a $(1,\delta)$-broom
with handle $w$, and $G[\{w,u_j,v_j\}\cup Q_j]$ is a $(2,\delta)$-broom with handle $w$; and taking the first for
$\delta$ choices of $j\in J_3$, and the second for the remaining $\delta$ choices of $j$, and taking their union,
we find that $G$ contains $T(\delta)$,
a contradiction. This proves (1).

\bigskip
By (1), for each $i\in \{1\ll n\}$ there is a subset $J_i\subseteq \{i+1\ll n\}$ with $|J_i|< 3r$, such that
for all $a,b,c$, if $i$ is strong to $(a,b,c)$ then one of $a,b,c$ is in $J_i$. Let $D$ be the digraph with vertex set
$\{1\ll n\}$ in which for $i<j$, $j$ is adjacent from $i$ if $j\in J_i$. Thus $D$ has maximum outdegree less than $3r$, and 
since $D$ is acyclic, 
by \ref{digraph} the graph underlying $D$ has chromatic number at most $3r$. It follows that there is a subset $I\subseteq \{1\ll n\}$
such that 
$$3r\chi(\bigcup_{i\in I}B_i\setminus \Pi)\ge \chi(\bigcup_{1\le i\le n}B_i\setminus \Pi),$$
with the property that for all $i,a,b,c\in I$, $i$ is not strong to $(a,b,c)$.
Let $W= \bigcup_{i\in I}B_i\setminus \Pi$. 
\\
\\
(2) {\em $\chi(G[W])\le s^2\delta\zeta\beta\tau^2$.}
\\
\\
Let $D$ be the digraph obtained from $G[W]$ by deleting all edges $uv$ such that $\{u,v\}$ is a subset of some $B_i(i\in I)$,
and directing every remaining edge such that the head of every edge is later than its tail.
Let $X$ be the set of all $v\in W$, where $v\in B_i$ say, such that for some $j>i$,
there are $s\delta$ neighbours of $v$ in $B_j$, pairwise nonadjacent.
Suppose that $D[X]$ has a directed path with $s+1$ vertices,  say $v_0\cc v_{s}$ in order. Let $v_i\in B_{j_i}$ where $j_i\in I$;
then $j_0\le \cdots\le j_{s}$. For $0\le i\le s$, since $v_{i}\in X$, there exist $b_i\in I$ with $b_i>j_i$
and a stable subset $P_{i}$
of $B_{b_i}$ with $|P_{i}| = s\delta$, all adjacent to $v_{i}$. For $1\le i\le s$, let $Q_i$ be the set of vertices in $P_s$
that are adjacent to $v_i$ and not to $v_{i-1}$. Since $i-1$ is not strong to $(j_{i}, b_{i-1}, b_{i})$, it follows that
$|Q_i|<\delta$ for $1\le i\le s$, and so there is a vertex $v\in P_s$ that belongs to none of the sets $Q_1\ll Q_s$.
Consequently $v$ is adjacent to all of $v_0,v_1\ll v_s$, and hence has a neighbour in $B_j$ for $s+1$ different values of $j$,
contrary to the hypothesis. This proves that $D[X]$ has no directed path with $s+1$ vertices, and hence the graph
underlying $D[X]$ has chromatic number
at most $s$, from the Gallai-Roy theorem~\cite{gallai,roy}. 

For each $v\in W\setminus X$, there are fewer than $s$ values of $j\in I$ such that $B_j\setminus \Pi$ 
contains some vertex adjacent from $v$ in $H$,
by hypothesis. Moreover, since $v\notin X$, if $v\in B_i\setminus \Pi$ and some vertex in $B_j\setminus \Pi$ is 
adjacent from $v$ in $D$, then
$j>i$, and there do not exist $s\delta$ vertices in $B_j\setminus \Pi$, pairwise nonadjacent in $G$ and all 
adjacent to $v$ in $G$. Since
the set of neighbours of $v$ is $\tau$-colourable, it follows that fewer than $s\delta\tau$ vertices in $B_j\setminus \Pi$ are 
adjacent in $D$
from $v$; and so $H[W\setminus X]$ has maximum outdegree less than $(s-1)s\delta\tau$. Thus by \ref{digraph} the graph underlying $D[W\setminus X]$
has chromatic number at most $(s-1)s\delta\tau$ (since $D$ is acyclic). 

Consequently the graph underlying $D$ has chromatic number at most $s+ (s-1)s\delta\tau\le s^2\delta\tau$. Since for each $i\in I$,
$G[B_i\setminus \Pi]$ has chromatic number at most $\beta\zeta\tau$, 
it follows that $\chi(G[W])\le s^2\beta\delta\zeta\tau^2$. This proves (2).

\bigskip

From the choice of $I$, 
$3r\chi(W)\ge \chi(\bigcup_{1\le i\le n}B_i\setminus \Pi),$
and so from (2),
$$\chi(\bigcup_{1\le i\le n}B_i\setminus \Pi) \le 3rs\beta\delta\zeta\tau^2.$$ 
This proves \ref{manyB}.~\bbox

We deduce, finally:
\begin{thm}\label{mainthmagain}
There exists $c$ such that if $G$ satisfies {\rm (i)--(v)} then $\chi(G)\le c$.
\end{thm}
\Proof
Let $\eta = \alpha+ 2(\delta+1)^3(\epsilon+1)^2$ and $\zeta = \eta+\delta$.
Let $\phi$ and $s$ satisfy \ref{nested}. Let $r$ be as in \ref{manyB}.
Let $c_2=3rs\beta\delta\zeta\tau^2$, let $c_1=c_2+\delta\tau^2$, and let $c=\phi(c_1)$.

Let $G$ satisfy (i)--(v), and suppose that $\chi(G)>c$. By \ref{nested},
there is a $3$-cleaned $(\zeta,\eta)$-template array $\mathcal{T}$
in $G$, with sequence $(Y_i,H_i)\;(1\le i\le n)$, such that $\chi(U(\mathcal{T}))>c_1$, and
such that for each vertex $v\in U(\mathcal{T})$, there are fewer than
$s$ values of $i\in \{1\ll n\}$ such that some neighbour of $v$ in $U(\mathcal{T})$ has a neighbour in $H_i$.
By \ref{privatize}, there
is a 2-cleaned (in fact 3-cleaned) $(\zeta,\eta)$-template array $\mathcal{T}'$, with sequence $(Y_i,H_i')\;(1\le i\le n)$
and a privatization $\Pi$ for $\mathcal{T}'$
such that
\begin{itemize}
\item $H_i'\subseteq H_i$ for $1\le i\le n$;
\item $U(\mathcal{T}')\subseteq U(\mathcal{T})$; and
\item $\chi(U(\mathcal{T}')\setminus \Pi)\ge \chi(U(\mathcal{T}))-\delta\tau^2>c_2$.
\end{itemize}
It follows that for each vertex $v\in U(\mathcal{T}')$, there are fewer than
$s$ values of $i\in \{1\ll n\}$ such that some neighbour of $v$ in $U(\mathcal{T}')$ has a neighbour in $H_i'$.
By \ref{manyB} applied to $\mathcal{T}'$, 
$$\chi(U(\mathcal{T}')\setminus \Pi)\le 3rs\beta\delta\zeta\tau^2,$$
a contradiction. This proves \ref{mainthmagain}.~\bbox

Consequently, this completes the proof of \ref{mainthm2}, and hence of \ref{mainthm1}.


\begin{thebibliography}{99}
\bibitem{cuttingedge} M. Chudnovsky, A. Scott and P. Seymour,
``Induced subgraphs of graphs with large chromatic number. XII. Distant stars'',
submitted for publication, {\tt arXiv:1711.08612}.
\bibitem{gallai}  T. Gallai, ``On directed graphs and circuits'', in {\em Theory of Graphs} (Proc.
Colloquium Tihany, 1966), Academic Press, 1968, 115--118.
\bibitem{gyarfas} A. Gy\'arf\'as, ``On Ramsey covering-numbers'',
{\em Coll. Math. Soc. J\'anos Bolyai}, in {\em Infinite and Finite Sets},
North Holland/American Elsevier, New York (1975), 10.
\bibitem{gyarfas2} A. Gy\'arf\'as, E. Szemer\'edi and Zs. Tuza, ``Induced subtrees in graphs
of large chromatic number'', {\em Discrete Math.} 30 (1980), 235--344.
\bibitem{kierstead1} H. A. Kierstead and S. G. Penrice, ``Radius two trees specify $\chi$-bounded classes'',
{\em J. Graph Theory} 18 (1994), 119--129.
\bibitem{rodl} H. Kierstead and V. R\"odl,   ``Applications of hypergraph coloring to coloring graphs not inducing
certain trees'', {\em Discrete Math.} 150 (1996), 187--193.
\bibitem{kierstead2} H. A. Kierstead and  Y. Zhu, ``Radius three trees in graphs with large chromatic number'',
{\em SIAM J. Disc. Math.} 17 (2004), 571--581.
\bibitem{roy} B. Roy,  ``Nombre chromatique et plus longs chemins d'un graphe'' ,
{\em Rev. Fran\c{c}aise Informat. Recherche Op\'erationnelle} 1 (1967), 129--132.
\bibitem{scott}
A. D. Scott, ``Induced trees in graphs of large chromatic number'',
{\em J. Graph Theory} 24 (1997), 297--311.
\bibitem{sumner} D. P. Sumner, ``Subtrees of a graph and chromatic number'', in
{\em The Theory and Applications of Graphs}, (G. Chartrand, ed.),
John Wiley \& Sons, New York (1981), 557--576.
\end{thebibliography}
\end{document}